\documentclass[12 pt]{amsart}

\usepackage{graphicx}
\usepackage{amsmath, amssymb, tikz, pgfplots, stmaryrd}
\usetikzlibrary{calc,plotmarks}
\usepgfplotslibrary{groupplots}
\usepackage{amsopn}
\usepackage[caption=false]{subfig}
\usepackage{todonotes}

\theoremstyle{plain}
\newtheorem{theorem}{Theorem}
\theoremstyle{remark}
\newtheorem{definition}[]{Definition}
\newtheorem{algorithm}[]{Algorithm}

\newcommand{\gradx}{\nabla_{\!x}}
  
\newcommand{\an}[1]{\alpha^{{(#1)}}}
\newcommand{\mn}[1]{\mu^{{(#1)}}}
\newcommand{\RRR}{\mathbb{R}}
\newcommand{\cO}{\mathcal{O}}
\newcommand{\cT}{\mathcal{T}}

\newcommand{\om}{\varOmega}
\newcommand{\vphi}{\varphi}
\newcommand{\vPhi}{\varPhi}
\newcommand{\that}{\hat{t}}
\newcommand{\vpatch}{{\omega_V}}
\newcommand{\vp}{{\vpatch}}
\newcommand{\tmax}{t_{\mathrm{max}}}
\newcommand{\pp}{{\scriptstyle{P}}}
\newcommand{\tT}{{\scriptstyle{T}}}
\newcommand{\pd}{\partial}
\newcommand{\divx}{{\mathrm{div}}_x}
\newcommand{\divxt}{{\mathrm{div}}_{x,t}}
\newcommand{\divxthat}{{\mathrm{div}}_{x,\that}}
\newcommand{\intd}{\mathrm{d}}
\newcommand{\tA}{\tilde A}
\newcommand{\tM}{\tilde M}
\newcommand{\Mnorm}[2]{\|#1\|_{M(#2)}}

\newcommand{\M}{\mathbb{M}}

\newenvironment{sark_tableau}[1][1]{%
  \setlength{\arrayrulewidth}{.8pt}%
  \tabular%
}{%
  \endtabular
}

\begin{document}

\title{Structure aware Runge-Kutta time stepping for spacetime tents}

\thanks{This work was funded in part by NSF grant DMS-1912779 
    and by Austrian Science Fund (FWF) grant F65:
    {\em Taming Complexity in Partial Differential Equations}.}


\author[J.~Gopalakrishnan]{Jay Gopalakrishnan}
\address{Portland State University, PO Box 751, Portland OR 97207,USA }
\email{gjay@pdx.edu}

\author[J.~Sch\"oberl]{Joachim Sch{\"oberl}}
\address{Technische Universit{\"a}t Wien, Wiedner Hauptstra\ss e 8-10, 1040 Wien, Austria}
\email{joachim.schoeberl@tuwien.ac.at}

\author[C.~Wintersteiger]{Christoph Wintersteiger}
\address{Technische Universit{\"a}t Wien, Wiedner Hauptstra\ss e 8-10, 1040 Wien, Austria}
\email{christoph.wintersteiger@tuwien.ac.at}

\begin{abstract}
    We introduce a new class of Runge-Kutta type methods suitable for
  time stepping to propagate hyperbolic solutions within tent-shaped
  spacetime regions.  Unlike standard Runge-Kutta methods, the new
  methods yield expected convergence properties when standard high
  order spatial (discontinuous Galerkin) discretizations are used.
  After presenting a  derivation of nonstandard order conditions
  for these methods, we show numerical examples of nonlinear
  hyperbolic systems to demonstrate the optimal convergence rates.  We
  also report on the discrete stability properties of these methods
  applied to linear hyperbolic equations.


\end{abstract}

\maketitle

\section{Introduction}
\label{sec:introduction}

For simulating wave phenomena, the state of the art relies heavily on
efficient and accurate numerical solution techniques for hyperbolic
systems. This paper is concerned with those solution techniques that
proceed by subdividing the spacetime into tent-shaped subregions
satisfying a causality condition.  Just as light cones are often used
to delineate what is causally possible and impossible in the
spacetime, tent-shaped spacetime regions are natural to impose
causality when numerically solving hyperbolic equations.  By
constraining the height of the tent pole, erected vertically in an
increasing time direction, one can ensure that the tent encloses the
domain of dependence of all its points. This constraint on the tent
pole height is a causality condition that a numerical scheme using
such tents should satisfy.  The spacetime subdivision into tents may
be unstructured, thus allowing such schemes to advance in time by
different amounts at different spatial locations, i.e., local time
stepping can be naturally built in while subdividing the spacetime
into tents.

The main contribution of this paper is a new explicit Runge-Kutta type
time stepping scheme for solving hyperbolic systems within a spacetime
tent. Standard time stepping methods cannot be directly applied on
tents, since tents are generally not a tensor product of a spatial
domain with a time interval.  A non-tensor product spacetime tent can
be mapped to a tensor product spacetime cylinder using a degenerate
Duffy-like transformation. This is the basis of the Mapped Tent
Pitching (MTP) schemes that we introduced previously in~\cite{mtp}. As
shown there, a spacetime Piola map can be used to pull back the
hyperbolic system from the tent to the spacetime cylinder.  Being a
tensor product domain, the spacetime cylinder, admits the use of
standard explicit time stepping schemes, like the classical RK4
scheme. However, as we shall show here, expected convergence rates are
{\em not} observed when such standard explicit schemes are used. The
cause of this problem can be traced back to the degeneracy of the
map. After illustrating this problem, we shall introduce a new
Structure Aware Runge-Kutta (SARK) scheme, which overcomes this problem.

We first reported the above-mentioned order reduction
in~\cite{mtp_sat}, where a fix was proposed for {\em linear}
hyperbolic systems, called the Structure Aware Taylor (SAT) scheme. In
contrast, the new SARK schemes of this paper are applicable to {\em
  both linear and nonlinear} hyperbolic systems.

Prior work on tent-based methods spans both the computational
engineering literature~\cite{Haber2006,MILLER2008194} and numerical
analysis literature~\cite{Falk1999,Monk2005}. These works have clearly
articulated the promise of tent-based schemes, including local time
stepping, even with higher order spatio-temporal discretizations, and
the opportunities to utilize concurrency.  Recent advances include
tent-based Trefftz methods~\cite{PerugSchobStock20} and the use of
asynchronous SDG (spacetime discontinuous Galerkin) methods to new
engineering applications~\cite{AbediHaber18}.

To place the present
contribution in the perspective of these existing works, a few words
regarding our focus on explicit time stepping are in order. The ratio of
memory movements to flops is very low for explicit schemes, making them
highly suitable  for the newly emerging many-core processors. However,
before the introduction of MTP schemes in~\cite{mtp}, it was not clear
that such advantages of explicit time stepping could be brought to any
tent-based method.  Now that we have an algorithmic avenue to perform
explicit time stepping within tent-based schemes, we turn to the
study of accuracy and convergence orders. Having encountered the
unexpected roadblock of the above-mentioned convergence order
reduction, we have been focusing on developing time stepping techniques to
overcome it. This paper is an outgrowth of these studies.

In the next section, we quickly review the construction of MTP
schemes, showing how the main system of ordinary differential
equations (ODEs) that is the subject of this paper arises. In
Section~\ref{sec:prob}, we show why we should not use standard
Runge-Kutta schemes for solving the ODE system. Then, in
Section~\ref{sec:rk-type}, we propose our new SARK schemes for solving
the ODE system. Section~\ref{sec:order-cond} derives order conditions
for these schemes. In
Section~\ref{sec:discrete-stability}, we  study  the discrete 
stability of SARK schemes.
Section~\ref{sec:numerical} reports on the good 
performance of the new schemes when applied to some
standard nonlinear hyperbolic systems.

\section{Construction of mapped tent pitching schemes}
\label{sec:mapped-tent-pitching}

In this section we  give a brief overview of MTP
schemes. A fuller exposition can be found in \cite{mtp}.
Let $\cT$ be a simplicial conforming spatial mesh of
a bounded spatial domain $\om_0 \subset \RRR^N.$
Spacetime tents are
built atop this spatial mesh in a sequence of steps. At the $i$th
step, a tent $K_i$ is added. It takes the form
\begin{equation}
  \label{eq:def_Ki}
  K_i := \{(x,t) : x\in\vpatch, \tau_{i-1}(x)\le t\le\tau_{i}\}
\end{equation}
where $\vpatch$ is the vertex patch made up of all elements in $\cT$
connected to a mesh vertex $V$. In~\eqref{eq:def_Ki},  the
function $\tau_i(x)$ is a continuous 
function of the spatial coordinate $x$ that is piecewise linear with
respect to $\cT$.
The graph of 
$\tau_i$ represents the advancing spacetime front at the $i$th step,
so the tent $K_i$ may be thought of as
the spacetime domain between these advancing fronts.
Within $K_i$, the distance the
central vertex $V$ can advance in time is restricted by the {\em  causality
constraint}
\begin{equation}
  \label{eq:causality_constraint}
  |\nabla \tau_i| < \frac{1}{c_{\max}}\,,
\end{equation}
where $c_{\max}$ is an upper bound for the local wavespeed
on $\vpatch$  of a hyperbolic system under consideration.

The hyperbolic systems we have in mind are general systems with $L$ unknowns in $N$
spatial dimensions, posed on the spacetime cylinder 
$\om := \om_0 \times (0,\tmax)$ for some
final time $\tmax$. For sufficiently regular
functions $g:\om\times\RRR^L \rightarrow \RRR^{L}$ and
$f:\om\times\RRR^L \rightarrow \RRR^{L\times N}$, the hyperbolic
problem is to find $u:\om \rightarrow \RRR^{L}$ such that
\begin{equation}
  \label{eq:conslaw}
  \pd_tg(x,t,u) + \divx f(x,t,u) = 0
\end{equation}
where $\pd_t = \pd\slash{\pd t}$ denotes the time derivative and $\divx$
denotes the row-wise divergence operator.
Hyperbolicity of~\eqref{eq:conslaw} implies that there is   a set of
eigenvalues---whose magnitudes give wavespeeds---for each direction vector and each
point $x,t,u$. Let $c(x,t,u)$ denote the maximum of these
wavespeeds over all directions.  The quantity $c_{\max}$ in
\eqref{eq:causality_constraint} can be taken to be any upper bound for
these $c(x, t, u)$ for $(x, t) \in K_i$.

MTP schemes proceed by mapping each of the tents arising above to a
spacetime cylinder. 
To define the mapping, we consider a general tent $K$ over any given vertex
patch $\vpatch$, defined by 
\begin{equation}
  \label{eq:def_tent}
  K := \{(x,t) : x\in\vpatch, \vphi_b(x)\le t\le\vphi_t(x)\}.
\end{equation}
The functions $\vphi_b$ and $\vphi_t$ are continuous functions
that are piecewise linear on the vertex patch and may be identified as 
the bottom and top advancing
fronts restricted to the vertex patch $\vpatch$. To map the tent $K$
to the spacetime cylinder $\hat K := \vpatch\times(0,1)$, we define
the transformation $\vPhi : \hat K \rightarrow K$ by 
$  \vPhi(x,\that) := (x,\vphi(x,\that)),$
where 
\begin{equation}
  \label{eq:def_phi}
  \vphi(x,\that) := (1-\that)\vphi_b(x) + \that\vphi_t(x).
\end{equation}
Defining 
$F:\om\times\RRR^L \rightarrow \RRR^{L\times (N+1)}$ by
\begin{equation*}
  F(x,t,u) := [f(x,t,u),g(x,t,u)] \in \RRR^{L\times (N+1)},
\end{equation*}
we may write (\ref{eq:conslaw}) as
\begin{equation*}
  \divxt F(x,t,u) = 0.
\end{equation*}
The spacetime divergence $\divxt$ is a row-wise operator which applies
the spatial derivative to the first $N$ components and the temporal
derivative to the last component. 
The well-known Piola transformation of $F$, defined by
$\hat F = \left(\det\vPhi'\right) \left(F\circ\vPhi\right)
\left(\vPhi'\right)^{-T}$ can be simplified after calculating the
derivatives of $\varPhi$ to
\begin{equation*}
  \hat F
  = \left(F\circ\vPhi\right) 
  \begin{bmatrix}
    \delta I & -\nabla\vphi \\ 0 & 1
  \end{bmatrix}
\end{equation*}
where $\delta(x) = \vphi_t(x)-\vphi_b(x)$ 
and $I\in\RRR^{N\times N}$ is the identity matrix.  By the properties
of the Piola map, we then immediately have
\begin{equation}
  \label{eq:mapped_conslaw_st}
  \divxthat \hat F(x,\that,\hat u) = 0,
\end{equation}
with $\hat u = u\circ\vPhi$ and the spacetime divergence $\divxthat$
on $\hat K$. Finally, as described in~\cite{mtp}, writing $\hat F$ in
terms of $f$ and $g$, we find that~(\ref{eq:mapped_conslaw_st}) is
equivalent to
\begin{equation}
  \label{eq:mapped_conslaw_1}
  \pd_{\that}\left(g(x,\that,\hat u) - f(x,\that,\hat u) \nabla\vphi(\that)\right) 
  + \divx\left(\delta(x)f(x,\that,\hat u)\right) = 0,
\end{equation}
which is again a conservation law.

For readability, we omit the spatial
variable $x$ and pseudo-time $\that$ from the arguments of functions
in~(\ref{eq:mapped_conslaw_1}) and simply write
\begin{equation}
  \label{eq:mapped_conslaw}
  \pd_{\that}\left(g(\hat u) - f(\hat u) \nabla\vphi\right) 
  + \divx\left(\delta f(\hat u)\right) = 0,
\end{equation}
which describes the evolution of $\hat u$ along pseudo-time from
$\that=0$ to $\that=1$.
Since
\begin{equation*}
  \vphi(x,\that) = (1-\that)\vphi_b(x) + \that\vphi_t(x) = \vphi_b(x) + \that \delta(x),
\end{equation*}
we may split $g(\hat u) - f(\hat u) \nabla\vphi$ into
parts with and without explicit dependence on pseudo-time,
allowing us to rewrite~\eqref{eq:mapped_conslaw}
as
\begin{equation*}
  \pd_{\that}\left(\left(g(\hat u) - f(\hat u) \nabla\vphi_b\right)
    - \that f(\hat u) \nabla \delta\right) 
  + \divx\left(\delta f(\hat u)\right) = 0.
\end{equation*}




This equation is the starting point for our spatial discretization. 
We use a 
discontinuous Galerkin method based on
\[
  V_h = \{ v: v|_K \text{ is a polynomial of degree } \le p
  \text{ on all spatial elements } T \in \cT\}.
\]
When restricted to the vertex patch $\vpatch$ we obtain
$V_h(\vpatch) = \{ v|_{\vpatch}: v \in V_h\}$.  Multiplying
(\ref{eq:mapped_conslaw}) by a test function $v_h \in V_h$ and integrating by
parts over the patch $\vpatch$, we obtain
\begin{equation}
  \label{eq:var_prob}
  \begin{aligned}
    \int_\vpatch\frac{ d}{ d\hat t}
    \bigg(g(\hat u) &  - f(\hat u) \nabla\vphi\bigg)
  \cdot v_h
  \\
  & =  \sum_{T\subset\vpatch}\int_T \delta f(\hat u) : \nabla v_h
  - \sum_{F\subset\vpatch}\int_F \delta f_n(\hat u^+,\hat u^-)
  \cdot \llbracket v_h \rrbracket,        
  \end{aligned}
\end{equation}
for all $v_h\in V_h$ and all $\that\in[0,1]$. Here and throughout,
every facet $F$ is assigned a unit normal, simply denoted by $n$,
whose direction is arbitrarily fixed, except when
$F \subset \partial\om$, in which case it points outward.  The
traces $\hat u^+$ and $\hat u^-$ of $\hat u$ from either side are
defined by 
\begin{equation*}
  \hat u^+ := \lim_{s\rightarrow 0^+}\hat u(x+sn) \quad \text{and} \quad 
  \hat u^- := \lim_{s\rightarrow 0^+}\hat u(x-sn).
\end{equation*}
In~\eqref{eq:var_prob}, we also used a numerical flux $f_n$ on each facet
$F$ (that takes values in $\RRR^L$ depending on values
$\hat u^+,\hat u^-$ from either side) and the jump
$\llbracket \hat v_h \rrbracket := \hat v_h^+ - \hat v_h^-$. In these
definitions, whenever $\hat{u}^+$ falls outside $\om$, it is
prescribed using some given boundary conditions.

Let $m = \dim V_h(\vpatch)$ and let $\psi_i$, $i=1, \ldots, m$ denote
any standard local basis for $V_h(\vpatch)$.
Introducing $U: [0,1] \to \RRR^m$,
consider the basis expansion 
\begin{equation}
  \label{eq:basis_exp}
  \hat u (x, \that ) =
  \sum_{i=1}^m U_i(\that) \psi_i(x).  
\end{equation}
Equation~\eqref{eq:var_prob} leads to an ODE system for $U(\that)$ as
follows.  Define (possibly nonlinear) operators
$M_0: \RRR^m \to V_h(\vp)$, $M_1: \RRR^m \to V_h(\vp)$, and
$A : \RRR^m \to V_h(\vp)$ by
\begin{subequations} \label{eq:MA}
  \begin{align}
    \label{eq:Mo}
  \int_\vp M_0(U) v_h
  & = \int_\vp
  \left(g(\hat u)  - f(\hat u) \gradx\vphi_b\right)
    \cdot v_h
  \\ \label{eq:M1}
  \int_\vp M_1(U) v_h
  & = \int_\vp f(\hat u) \gradx \delta\; \cdot \,v_h
  \\ 
  \int_\vp A(U) v_h
  & =  \sum_{T\subset\vpatch}\int_T \delta f(\hat u) : \nabla v_h
  - \sum_{F\subset\vpatch}\int_F \delta \, f_n(\hat u^+,\hat u^-)
    \cdot \llbracket v_h \rrbracket,    
\end{align}
\end{subequations}
for all $v_h \in V_h(\vp)$, where $\hat u$ in all terms on the right
hand sides above is to be expanded in terms of $U$
using~\eqref{eq:basis_exp}.

With these notations, problem~\eqref{eq:var_prob} becomes the problem of
finding $U: [0,1] \to \RRR^m$ satisfying
\begin{equation}
  \label{eq:discr_ode}
  \frac{d}{d\that} M(\that, U(\that)) = A (U(\that)), 
\end{equation}
given some $U(0) = U_0 \in \RRR^m.$  Here  $M : [0, 1] \times \RRR^m \to V_h(\vpatch)$ is defined by
\begin{align}
  M(\that, W) = M_0(W) - \that M_1(W), \qquad 0 \le \hat t \le 1, \quad
  W \in \RRR^m. \label{eq:Mt}
\end{align}

\section{Difficulty with standard time stepping}
\label{sec:prob}

In this section, we describe the problem we must overcome,
thus setting the stage for the new schemes proposed in
Section~\ref{sec:rk-type}. The problem is that
standard Runge-Kutta methods when applied to the tent
system~\eqref{eq:discr_ode}---after a standard reformulation---do not
give expected orders of convergence.

A standard approach to numerically solve~\eqref{eq:discr_ode} proceeds
by introducing a new variable $Y(\that) := M(\that, U(\that)).$ Then,
using the inverse mapping $M^{-1}(\hat t, \cdot)$, the primary
variable is $U(\hat t) = M^{-1}(t, Y(\hat t))$. Substituting this
expression for $U$ on the right hand side of~\eqref{eq:discr_ode}, we
can bring~\eqref{eq:discr_ode} to the standard form
\begin{equation}
  \label{eq:discr_ode_subs}
  \frac{d}{d\that}Y(\that) - A(M^{-1}(\hat t, Y(\that))) = 0.
\end{equation}
Standard ODE solvers, such as the classical explicit RK4 method, may
now be directly applied to~\eqref{eq:discr_ode_subs}.  Unfortunately, this leads
to reduced convergence order, as we shall now see.

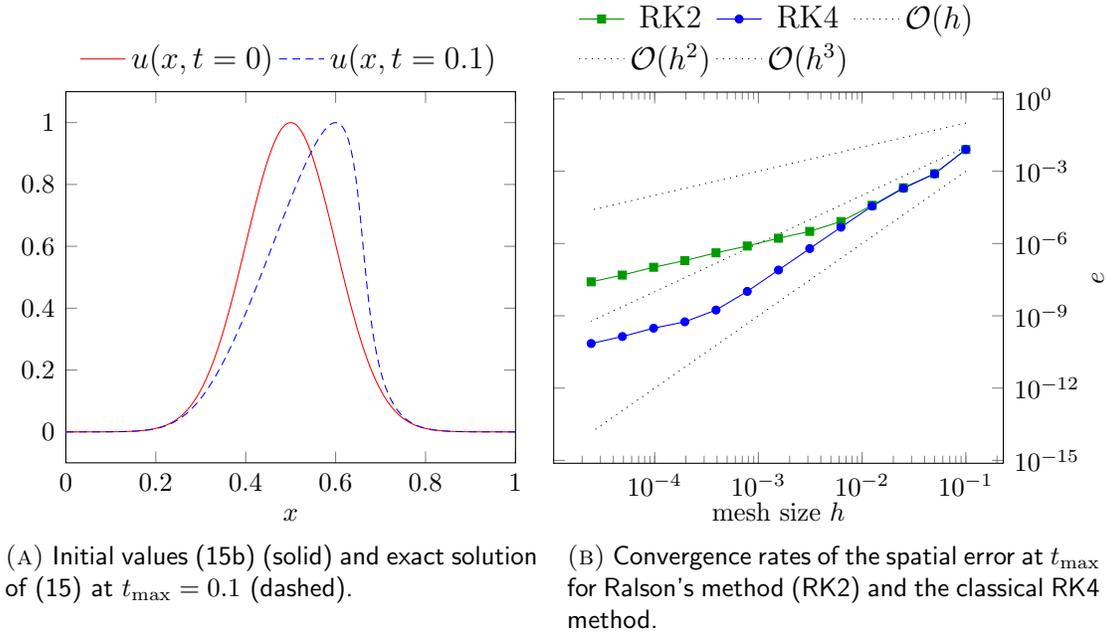
\begin{figure}
  \captionsetup{width=.45\linewidth}
  \centering
  \subfloat[Initial values (\ref{eq:burgers_example_u0}) (solid)
  and exact solution of (\ref{eq:burgers1d_example}) at $t_{\max}=0.1$ (dashed).]{
    \label{fig:burgers1d_inital_exact}
    \begin{tikzpicture}
      \begin{axis}[%
        width=0.48\textwidth,%
        tick label style={font=\footnotesize},%
        xmax=1, xmin=0,%
        legend style = {at={(0.5,1.02)}, anchor=south, draw=none},%
        legend columns = 2,
        xlabel={\footnotesize $x$},
        xlabel style = {at={(0.5,0.01)} ,anchor=north},
        ]%
        \addplot[color=red] table[x=x,y=u0] {burgers1d_sol_t0_1.txt};
        \addlegendentry{$u(x,t=0)$};
        
        \addplot[color=blue,style=densely dashed] table[x=x,y=u0_1] {burgers1d_sol_t0_1.txt};
        \addlegendentry{$u(x,t=0.1)$};
      \end{axis}
    \end{tikzpicture}
  }
  \subfloat[Convergence rates of the spatial error at $t_{\max}$ for Ralson's method (RK2) and the classical RK4 method. ]{
    \label{fig:burgers1d_rk_rates}
    \begin{tikzpicture}
      \tikzset{mark options={mark size=1.5}}%
      \begin{loglogaxis}[%
        width=0.48\textwidth,%
        tick label style={font=\footnotesize},%
        yticklabel pos=right,
        legend style = {at={(0.5,1.02)}, anchor=south, draw=none},%
        legend columns = 3,
        xlabel={\footnotesize mesh size $h$},%
        xlabel style = {at={(0.5,0.03)} ,anchor=north},%
        ylabel={\footnotesize $e$},%
        ylabel style = {at={(1.38,0.5)} ,anchor=north},%
        ]%
        \addplot[color=green!60!black,mark=square*] table[x=h,y=l2e_RK2_ralston] {table_burgers1d_conv_rk.txt};
        \addlegendentry{RK2};
        \addplot[color=blue,mark=*] table[x=h,y=l2e_RK4] {table_burgers1d_conv_rk.txt};
        \addlegendentry{RK4};
        
        \addplot[dotted] table[x=h,y=h] {table_burgers1d_conv_rk.txt};
        \addlegendentry{$\cO(h)$};
        \addplot[dotted] table[x=h,y expr={x^2}] {table_burgers1d_conv_rk.txt};
        \addlegendentry{$\cO(h^2)$};
        \addplot[dotted] table[x=h,y expr={x^3}] {table_burgers1d_conv_rk.txt};
        \addlegendentry{$\cO(h^3)$};
    \end{loglogaxis}
  \end{tikzpicture}
}
\caption{Exact solution $u$ and convergence rates of the error $e$,
  defined in (\ref{eq:def_l2e}), for the example of the Burger's
  equation described in (\ref{eq:burgers1d_example}).}
\label{fig:burgers1d}
\end{figure}

Consider the example of the one-dimensional Burger's equation
\begin{subequations}
  \label{eq:burgers1d_example}
  \begin{equation}
    \pd_{t} u(x,t) + \pd_x u(x,t)^2 = 0, \quad \forall (x,t)\in[0,1]\times(0,t_{\max}],
  \end{equation}
  with initial values set by 
  \begin{equation}
    \label{eq:burgers_example_u0}
    u(x,0)=\exp{\big(-50\big(x-\tfrac{1}{2}\big)^2\big)}, \quad \forall x\in[0,1],
  \end{equation}
an inflow boundary condition an $x=0$, and an outflow boundary
condition at $x=1$. The final time $t_{\max} = 0.1$ is chosen such
that the exact solution is still a smooth function (before the onset
of shock).
Therefore no regularization or limiting
is expected to be essential to witness  high order
convergence. The exact solution at $t_{\max}$ shown in
Fig.~\ref{fig:burgers1d_inital_exact} is obtained by the method of
characteristics together with a Newton's method.  Let $u_h(x)$ be the
numerical solution of (\ref{eq:burgers1d_example}) at $t=t_{\max}$
computed with an explicit MTP scheme using a DG method of order
$p = 2$ for the spatial discretization. Thus one would   expect the
error
\begin{equation}
  \label{eq:def_l2e}
  e := \|u(\cdot,0.1)-u_h\|_{L^2([0,1])}  
\end{equation}
\end{subequations}
to go to zero at a rate of $\cO(h^3)$. However our observations in
Figure~\ref{fig:burgers1d_rk_rates} run counter to that expectation.
Figure~\ref{fig:burgers1d_rk_rates} reports the rates we observed when
two standard time stepping schemes were used to solve~(\ref{eq:discr_ode_subs}),
namely the two-stage (RK2) and the four-stage (RK4) explicit
Runge-Kutta time stepping schemes.  Although we see third order
convergence for the first few refinement steps, the rate eventually
drops to first order for both methods.
We shall return to this example in
Subsection~\ref{ssec:conv-rates-burg}
after developing a method without this convergence reduction.

\section{Structure aware Runge-Kutta type methods}
\label{sec:rk-type}

In this section, we develop specialized Runge-Kutta type schemes that
do not show the above mentioned order loss of classical Runge-Kutta
schemes. The appearance of $M^{-1}$ in~\eqref{eq:discr_ode_subs} is
the cause of the problem.  Although $M$ is linear in $\hat t$, its
inverse is not. Therefore a spacetime polynomial solution of a
hyperbolic problem is not generally preserved even by high order
discretizations of~\eqref{eq:discr_ode_subs}.

With this in mind, we  motivate the definition of the new scheme
by reformulating~\eqref{eq:discr_ode} in terms of two 
variables $Z(\hat t)$ and $Y(\hat t)$, defined by
\[
  Z(\hat t) = M_0(U(\hat t)), \qquad 
  Y(\hat t) = M(\hat t, U(\hat t)) = Z(\hat t) - \that M_1(U(\that)).
\]
Then~\eqref{eq:discr_ode} implies
\begin{equation}
  \label{eq:discr_ode_yz}
  Y' = A(U(\hat t)), \qquad Z' =
  A(U(\hat t))  + (\hat t M_1(U(\hat t)))',
\end{equation}
together with the initial conditions $ Y(0) = Z(0) = M_0(U_0).$
Here and throughout we use primes ($'$) to abbreviate $d / d\hat t$.
The key idea is to avoid the inversion of the time-dependent~$M$
at all $\hat t$, limiting the inversion to just that of
$M_0$. Assuming we can compute the time-independent inverse
$M_0^{-1}$, we define
\[
  \tA = A \circ M_0^{-1}, \qquad \tM_1 = M_1 \circ M_0^{-1}.
\]
Then,~\eqref{eq:discr_ode_yz} yields the following ODE system for $Y$ and $Z$ on
$0 < \hat t < 1$:
\begin{subequations}
  \label{eq:exactflow}
  \begin{align}
    \label{eq:exactflow-Z}
    Z' &  = \tA( Z(\hat t))  + (\hat t \tM_1(Z(\hat t)))',
    && Z(0) = Y_0,\\
    \label{eq:exactflow-Y}
    Y' &  = \tA( Z(\hat t)),
    && Y(0) = Y_0,
  \end{align}
\end{subequations}
where $Y_0 = M_0(U_0)$.

Integrating the equations of~\eqref{eq:exactflow} from $0$ to
$\tau$, we obtain
\begin{subequations}
  \begin{align}
    \label{eq:rktype_z}
    Z(\tau) & = Z(0) + \tau \tM_1(Z(\tau)) + \int_0^{\tau} \tA(Z(s)) \,\intd s,
    \\
    \label{eq:rktype_y}
    Y(\tau) & = Y(0)+
              \int_0^{\tau} \tA (Z(s))\, \intd s.
  \end{align}
\end{subequations}
The new scheme, defined below, may be thought of as motivated by
quadrature approximations to the integrals above. Note that we are
only interested in such quadratures that result in explicit schemes.
Moreover, we must also approximate $\tau \tM_1(Z(\tau))$ by an
extrapolation formula that uses prior values of $Z,$ in order to keep
the scheme explicit.

\begin{definition}
  Given an initial condition $Y_0$,  {\em an $s$-stage SARK method}
  for~\eqref{eq:exactflow} computes 
  \begin{subequations} \label{eq:rktype}
    \begin{align}
      Z_i & = Y_0 + \tau \sum_{j<i} d_{ij}\tM_1(Z_j)
            + \tau \sum_{j<i}a_{ij}\tA(Z_j), \qquad 1\le i \le s,
            \label{eq:rktype_Zi}
      \\
      Y_\tau & = Y_0 + \tau \sum_{i=1}^s b_i \tA(Z_i). \label{eq:rktype_int}
    \end{align}
  \end{subequations}
\end{definition}

This explicit method is determined by the coefficient matrices $b \in
\RRR^{s \times 1}$,
$\mathcal{A} \in \RRR^{s \times s}$,
and
$\mathcal{D} \in \RRR^{s \times s}$:
\begin{gather*}
  b = (b_1,\dots,b_s), \quad 
  \mathcal{A}  =
                \begin{pmatrix}
                  0 & & & \\
                  a_{21} & 0 &  &\\
                  \vdots & \ddots & 0 & \\
                  a_{s1} & \dots & a_{s,s-1} & 0\\
                \end{pmatrix}, \quad 
  \mathcal{D}  =
                \begin{pmatrix}
                  0 & & & \\
                  d_{21} & 0 &  &\\
                  \vdots & \ddots & 0 & \\
                  d_{s1} & \dots & d_{s,s-1} & 0\\
                \end{pmatrix}.
\end{gather*}
Hence we use
\begin{tabular}{c|c|c}
  $c$ & $\mathcal A$ & $\mathcal D$ \\ \hline &  $b$ &
\end{tabular} instead 
of the standard Butcher tableau 
\begin{tabular}{c|c}
  $c$ & $\mathcal A $ \\ \hline & $b$ 
\end{tabular}
to express our scheme. Here we restrict ourselves to schemes where
$c \in \RRR^s$ is set by the consistency condition
\[
  c_i= \sum_{j=1}^{i-1} a_{ij}.
\]
In the next section, we shall develop a theory to  choose appropriate values of
$a_{ij}, d_{ij},$ and $b_i$.  There, Subsection~\ref{sec:examples} contains
some specific examples of SARK scheme tableaus.

\section{Order conditions for the scheme}
\label{sec:order-cond}

Appropriate values of $a_{ij}, d_{ij},$ and $b_i$ can be found by
order conditions obtained by matching terms in the Taylor expansions of the
exact solution $Y(\tau)$ and the discrete solution $Y_\tau$.  To
derive these order conditions we follow the general methodology laid
out in~\cite{hairer2008}. For this, we need to first compute the
derivatives of the exact flow (in~\S\ref{ssec:deriv-exact}), then the
derivatives of the discrete flow (in~\S\ref{ssec:deriv-discr}),
followed by the formulation of resulting order conditions (in
\S\ref{ssec:conditions-order-3}).

\subsection{Derivatives of the exact solution}  \label{ssec:deriv-exact}

Continuing to use primes ($'$) for total derivatives with respect to a single
variable like $d/d \tau$, to ease the tedious calculations below, we
shall also employ the $n$th order Frechet derivative of a function
$g: D \subset \RRR^m \to V$, for some vector space $V$. It is denoted
by $g^{(n)}(z): \RRR^m \times \cdots \times \RRR^m \to V$ and defined
by the symmetric multilinear form
\[
  g^{(n)}(z)(v_1, \ldots, v_n) =
  \sum_{i_1, i_2, \ldots, i_n=1}^m
  \frac{ \pd^n g(z)}{ \pd x_{i_1} \cdots \pd x_{i_n}}
  [v_1]_{i_1} \ldots [v_n]_{i_n}
\]
for any $v_1, \ldots, v_n \in \RRR^m.$ Whenever $g$ and
$z: (0,1) \to \RRR^m$ are sufficiently smooth for the derivatives
below to exist continuously, we have the following formulae.
\begin{subequations}
  \label{eq:derivs_compos}
  \begin{align}
    \frac{d}{d \tau} g(z(\tau))
    & =  g^{(1)}(z(\tau))(z'(\tau)),
    \\
    \frac{d^2}{d \tau^2} g(z(\tau))
    & =  g^{(2)}(z(\tau))(z'(\tau), z'(\tau)) +
      g^{(1)}(z(\tau)(z''(\tau)),
    \\
    \frac{d^3}{d \tau^3} g(z(\tau))
    & =  g^{(3)}(z(\tau))(z'(\tau), z'(\tau), z'(\tau))
    \\ \nonumber 
    & + 3 g^{(2)}(z(\tau))(z'(\tau), z''(\tau))
      + g^{(1)}(z(\tau))(z'''(\tau)),
    \\ 
    \frac{d^4}{d \tau^4} g(z(\tau))
    & = g^{(4)}(z(\tau))(z'(\tau), z'(\tau), z'(\tau), z'(\tau))
    \\ \nonumber 
    & + 6 g^{(3)}(z(\tau))(z'(\tau), z'(\tau), z''(\tau))
      + 4 g^{(2)}(z(\tau))(z'(\tau), z'''(\tau))
    \\ \nonumber 
    & + 3 g^{(2)}(z(\tau)) (z''(\tau), z''(\tau)) + g^{(1)}(z(\tau))(z''''(\tau)).
  \end{align}
\end{subequations}
These formulae can be derived by repeated application of the chain
rule (or by applying the Fa{\'a} di Bruno formula). We will also
need to use
\begin{equation}
  \label{eq:Leibniz_tau}
  \frac{d^k}{d \tau^k}\big( \tau g(z(\tau)) \big)
  = \tau \frac{d^k}{d \tau^k} g(z(\tau)) + k
  \frac{d^{k-1}}{d \tau^{k-1}} g(z(\tau)),
\end{equation}
which is a simple consequence of the Leibniz rule.

We start by computing the derivatives of $Z(\tau)$ at
$\tau=0$. To express such derivatives 
concisely, we introduce the notation
\begin{gather*}
  \alpha = \tA(Z(0)),
  \qquad \an n (v_1, \ldots, v_n) = \tA^{(n)}(Z(0))(v_1, \ldots, v_n),
  \\
  \mu = \tM_1(Z(0)),
  \qquad \mn n (v_1, \ldots, v_n) = \tM^{(n)}(Z(0))(v_1, \ldots, v_n).
\end{gather*}
From~\eqref{eq:exactflow-Z}, it is immediate that $Z'(0) =
\tA (Z(0)) + \tM(Z(0))$. Thus, 
\begin{subequations}
  \label{eq:dZ}
  \begin{align}
    \label{eq:dZ_1}
    Z'(0) &= \alpha + \mu.
  \end{align}
  For the next derivative, we differentiate~\eqref{eq:exactflow-Z}
  twice to get $Z''(\tau) = (\tA(Z(\tau))' + (\tau \tM_1(Z(\tau)))''.$
  Calculating the latter using~\eqref{eq:Leibniz_tau}, simplifying
  using~\eqref{eq:derivs_compos}, and evaluating at  $\tau=0$, we obtain 
  \begin{equation}
    \label{eq:dZ_2}
    Z''(0) = (\an 1 + 2 \mn 1)(\alpha + \mu).
  \end{equation}
  By the same procedure, starting with
  $Z'''(\tau) = (\tA(Z(\tau))'' + (\tau \tM_1(Z(\tau)))'''$
  and using~\eqref{eq:Leibniz_tau}
  and~\eqref{eq:derivs_compos}, we also have
  \begin{equation}
    \label{eq:dZ_3}
    Z'''(0) = (\an 2 + 3\mn 2)(\alpha + \mu, \alpha + \mu) +
    (\an 1 + 3 \mn 1)
    \big(
    (\an 1 + 2 \mn 1 )(\alpha + \mu)
    \big).
  \end{equation}
\end{subequations}

Armed with~\eqref{eq:dZ}, we proceed to compute the derivatives of
$Y$. Obviously, \eqref{eq:exactflow-Y} implies 
\begin{subequations}  \label{eq:dY0}
  \begin{align}
    \label{eq:dY0_1}
    Y'(0) & = \tA(Z(0)) = \alpha.
\intertext{Differentiating \eqref{eq:exactflow-Y} again,
  using~\eqref{eq:derivs_compos}, and evaluating at $\tau=0$ using the
  previously computed derivatives of $Z$ in~\eqref{eq:dZ}, we also get}
            \label{eq:dY0_2}
    Y''(0) & =
             \an 1 (\alpha + \mu)
    \\
    \label{eq:dY0_3}
    Y'''(0) & = \an 2 (\alpha+\mu, \alpha+\mu) +
             \an 1
              \big(
              (\an 1 + 2 \mn 1 )(\alpha + \mu)
             \big).
  \end{align}
\end{subequations}

\subsection{Derivatives of the discrete flow} \label{ssec:deriv-discr}

The next task is to compute the coefficients of the Taylor expansion
of the function $Y_\tau$ defined in~\eqref{eq:rktype_int}. The
arguments $Z_i$ in~\eqref{eq:rktype_int} are also functions of $\tau$,
as given by~\eqref{eq:rktype_Zi}. Therefore, in what follows, we first
differentiate $Z_i \equiv Z_i(\tau)$ and then $Y_\tau$.

Obviously, $Z_i(0)$ and $Z(0)$ coincide, so we will focus on the first
and higher derivatives of $Z_i$ at $\tau=0$.  To this end, we
differentiate \eqref{eq:rktype_Zi} $k$ times to get
\[
  \frac{d^k Z_i}{d \tau^k} =  \sum_{j<i} 
  \bigg[ d_{ij}\frac{d^k}{ d\tau^k} (\tau \tM_1(Z_j(\tau)))+
  a_{ij}
  \frac{d^k}{ d \tau^k}(\tau \tA(Z_j(\tau)))\bigg].
\]
Using~\eqref{eq:Leibniz_tau} for $k=1,2, 3$,
then~\eqref{eq:derivs_compos}, and evaluating at $\tau=0$ we obtain
\begin{subequations} \label{eq:dZi}
  \begin{align}
    Z_i'(0) & = \sum_{j < i} d_{ij} \mu + a_{ij} \alpha
    \\
    Z_i''(0) & = 2 \sum_{j < i} \sum_{k < j} 
               \big(
               d_{ij}\mn 1 + a_{ij}\an 1 
               \big)
               (d_{jk} \mu + a_{jk} \alpha)
    \\
    %
    Z_i'''(0)
            & = 3 \sum_{j< i} \sum_{k < j} \sum_{l < j}
              \big(
              d_{ij} \mn 2 + a_{ij} \an 2
              \big)
              (
              d_{jk} \mu + a_{jk} \alpha,
              d_{jl} \mu + a_{jl} \alpha
              )
    \\ \nonumber
            & + 6 \sum_{j< i} \sum_{k < j} \sum_{l < k}
              \big(
              d_{ij} \mn 1 + a_{ij} \an 1 
              \big)
              \big(
              (d_{jk} \mn 1 + a_{jk} \an 1)
              (d_{kl} \mu + a_{lk}\alpha)
              \big).
  \end{align}
\end{subequations}

Next, we focus on $Y_\tau$. 
By~\eqref{eq:rktype_int},
\[
  \frac{d^k Y_\tau}{d \tau^k}
  = \sum_{i=1}^s b_i
  \frac{d^k }{d \tau^k}(\tA (Z_i(\tau))).
\]
Using~\eqref{eq:derivs_compos}, and evaluating the resulting terms at
$\tau = 0$ by means of~\eqref{eq:dZi}, we obtain
\begin{subequations}  \label{eq:dY10}
  \begin{align}
    \label{eq:dY10_1}
    Y_\tau'(0) & = \sum_{i=1}^s b_i \alpha,
    \\
    \label{eq:dY10_2}
    Y_\tau''(0) & = 2\sum_{i=1}^s \sum_{j < i} b_i
                  \an 1(d_{ij} \mu + a_{ij} \alpha),
    \\
    \label{eq:dY10_3}
    Y_\tau'''(0)
               & =
                 3 \sum_{i=1}^s\sum_{j< i} \sum_{k < i}
                 b_i
                 \an 2
                 (d_{ij} \mu  + a_{ij} \alpha,
                 d_{ik} \mu  + a_{ik} \alpha)
    \\ \nonumber 
               & +
                 6 \sum_{i=1}^s\sum_{j< i} \sum_{k < j}
                 b_i\an 1
                 \big(
                 (d_{ij}\mn 1 + a_{ij} \an 1 )
                 (d_{jk} \mu + a_{jk} \alpha)
                 \big).
  \end{align}
\end{subequations}

\subsection{Formulation of order conditions} \label{ssec:conditions-order-3}

To obtain a specific method, we find values for $a_{ij}, d_{ij}$ and
$b_i$ by matching the coefficients in the Taylor expansions of
$Y(\tau)$ and $Y_\tau$. Note that $Y_\tau(0) = Y_0 = Y(0)$, so the
$0$th order coefficients match. 

The next terms in the Taylor expansions will match if
$Y'(0)=Y_\tau'(0)$. For this it is sufficient that
\begin{align}
  \label{eq:cond_1}
  \sum_{i=1}^s b_i = 1.
\end{align}
because of~\eqref{eq:dY0_1} and~\eqref{eq:dY10_1}. To match the third
terms in the Taylor expansions, equating~\eqref{eq:dY0_2}
and~\eqref{eq:dY10_2},
\[
  \an 1 (\alpha) + \an 1 (\mu) =
  \sum_{i=1}^s \sum_{j < i} 2b_i d_{ij}
  \an 1 (\mu) + 2 b_ia_{ij} \an 1 (\alpha).
\]
Equating the coefficients of $\an 1 (\alpha)$ and $\an 1 (\mu)$, we
conclude that $Y''(0) = Y_\tau''(0)$ if
\begin{align}
  \label{eq:cond_2}
  2\sum_{i=1}^s \sum_{j<i}
  b_i d_{ij} = 1 \quad \text{and} \quad 2\sum_{i=1}^s\sum_{j<i} b_ia_{ij} = 1.
\end{align}

If one desires to further match the next higher order terms,
$Y_\tau'''(0) = Y'''(0)$, then the expressions in \eqref{eq:dY0_3}
and~\eqref{eq:dY10_3} must be equated, i.e.,
\begin{align*}
  \an 2 &(\alpha,\! \alpha)
  + 2 \an 2 (\alpha,\!\mu) + \an 2 (\mu,\! \mu) 
  \\
  & + \an 1 (\an 1 (\alpha)) + \an 1 (\an 1 (\mu))
  + 2 \an 1 (\mn 1 (\alpha) ) + 2 \an 1 (\mn 1 (\mu))
  \\
  & =  \sum_{i=1}^s \sum_{j< i} \sum_{k < i}
  \bigg[
  3 b_i d_{ij} d_{ik} \an 2 (\mu,\! \mu) + 6 b_i d_{ij} a_{ik} \an 2
  (\mu, \alpha) + 3 b_i a_{ij} a_{ik} \an 2 (\alpha, \alpha)
  \bigg]
  \\
   & + 6 \sum_{i=1}^s \sum_{j< i} \sum_{k < j}
     \bigg[
     b_i d_{ij}d_{jk} \an 1 (\mn 1 (\mu)) + b_i d_{ij}a_{jk} \an 1 (\mn 1 (\alpha))
  \\
   & \hspace{2.1cm} + 
     b_i a_{ij}d_{jk} \an 1 (\an 1 (\mu)) +
     b_i a_{ij}a_{jk} \an 1 (\an 1 (\alpha))   \bigg].            
\end{align*}
For this equality to hold, the following seven conditions are sufficient
as can be seen by equating the coefficients of
$\an2 (\alpha, \alpha),$ $\an 2 (\mu, \mu),$ $\an 2(\alpha, \mu),$ $\an 1
(\an 1 (\alpha)),$ $\an 1 (\an 1 (\mu)),$ $\an 1 (\mn 1 (\alpha) ),$ and
$ \an 1 (\mn 1 (\mu))$, respectively:
\begin{subequations}   \label{eq:cond_3}
\begin{align}
  3 \sum_{i=1}^s b_i \Big(\sum_{j<i} a_{ij}\Big)^2 = 1, \\
  3 \sum_{i=1}^s b_i \Big(\sum_{j<i} d_{ij}\Big)^2 = 1, \\
  3 \sum_{i=1}^s b_i \Big(\sum_{j<i} a_{ij}\Big)\Big(\sum_{j<i} d_{ij}\Big) = 1, \\
  6 \sum_{i=1}^s \sum_{j<i}\sum_{k<j} b_ia_{ij}a_{jk} = 1, \\
  6 \sum_{i=1}^s \sum_{j<i}\sum_{k<j} b_ia_{ij}d_{jk} = 1, \\
  3 \sum_{i=1}^s \sum_{j<i}\sum_{k<j} b_id_{ij}a_{jk} = 1, \\
  3 \sum_{i=1}^s \sum_{j<i}\sum_{k<j} b_id_{ij}d_{jk} = 1.
\end{align}
\end{subequations}
Thus, we have proved the following result, which summarizes our 
discussions on order conditions.

\begin{theorem}
  \label{thm:order-conditions}
  Whenever $\tA$ and $\tM$ are smooth enough for the derivatives below
  to exist continuously,
  \begin{enumerate}
  \item  the condition~\eqref{eq:cond_1} implies $ Y'(0) = Y_\tau'(0),$
  \item  the conditions of~\eqref{eq:cond_2} imply $Y''(0) =
    Y_\tau''(0),$ and 
  \item  the conditions of~\eqref{eq:cond_3} imply $ Y'''(0) = Y_\tau'''(0)$.
  \end{enumerate}
\end{theorem}

\subsection{Examples of methods up to third order}
\label{sec:examples}

Observe that the standard order conditions of Runge-Kutta methods
are a subset of the order
conditions derived in \S\ref{ssec:conditions-order-3}. Thus we
base our SARK methods on existing Runge-Kutta methods.
Below, we shall refer to an $s$-stage SARK method based on an existing
Runge-Kutta method called ``RKname'' as ``SARK($s$, RKname)''.

A second order two-stage SARK method can be derived
from a second order Runge-Kutta method once we find $d_{ij}$
satisfying the additional condition
\begin{align}
  \label{eq:sark2_coeff}
  2\sum_{i=1}^2 \sum_{j<i}
  b_i d_{ij} = 1 \quad \Leftrightarrow \quad b_2 d_{21} = \frac{1}{2},
\end{align}
which was introduced in (\ref{eq:cond_2}). For example, one may start
with the standard explicit midpoint rule and select $d_{21}=1/2$ to
satisfy~\eqref{eq:sark2_coeff}, thus arriving at the ``SARK(2, midpoint)''
method, listed first in 
Table~\ref{tab:sark2_examples}.
The table continues on to display
further such methods obtained from other well-known second order
Runge-Kutta schemes.

\begin{table}[b]
  \centering
  \subfloat[SARK(2, midpoint), based on the explicit midpoint rule]{
    \begin{sark_tableau}{c|cc|cc}
      0 & 0 & 0 & 0 & 0\\
      $\frac{1}{2}$ & $\frac{1}{2}$ & 0 & $\frac{1}{2}$ & 0\\
      \hline
      & 0 & 1 & &
    \end{sark_tableau}
    \label{tab:sark2_mprule}
  }
  \hspace{2em}
  \subfloat[SARK(2, Ralston), based on Ralston's second order method]{
    \begin{sark_tableau}{c|cc|cc}
      0 & 0 & 0 & 0 & 0\\
      $\frac{2}{3}$ & $\frac{2}{3}$ & 0 & $\frac{2}{3}$ & 0\\
      \hline
      & $\frac{1}{4}$ & $\frac{3}{4}$ & &
    \end{sark_tableau}
    \label{tab:sark2_ralston}
  }
  \hspace{2em}
  \subfloat[SARK(2, Heun), based on Heun's second order method]{
    \begin{sark_tableau}{c|cc|cc}
      0 & 0 & 0 & 0 & 0\\
      1 & 1 & 0 & 1 & 0\\
      \hline
      & $\frac{1}{2}$ & $\frac{1}{2}$ & &
    \end{sark_tableau}
    \label{tab:sark2_heun}
  }
  \caption{Two-stage SARK methods}
  \label{tab:sark2_examples}
\end{table}
\begin{table}
  \centering
  \subfloat[SARK(3, Kutta) method, based on Kutta's third order method]{
    \begin{sark_tableau}{c|rcc|rcc}
      0 & 0 & 0 & 0 & 0 & 0 & 0 \\
      $\frac{1}{2}$ & $\frac{1}{2}$ & 0 & 0 & $\frac{1}{2}$ & 0 & 0 \\
      1 & $-1$ & 2 & 0 & $-3$ & 4 & 0 \\
      \hline
      & $\frac{1}{6}$ & $\frac{2}{3}$ & $\frac{1}{6}$ &  &  &    
    \end{sark_tableau}
    \label{tab:sark3_kutta}
  }
  \hspace{2em}
  \subfloat[SARK(3, Heun) method, based on Heun's third order method]{
    \begin{sark_tableau}{c|ccc|rcc}
      0 & 0 & 0 & 0 & 0 & 0 & 0 \\
      $\frac{1}{3}$ & $\frac{1}{3}$ & 0 & 0 & $\frac{1}{3}$ & 0 & 0 \\
      $\frac{2}{3}$ & 0 & $\frac{2}{3}$ & 0 & $-\frac{2}{3}$ & $\frac{4}{3}$ & 0 \\
      \hline
      & $\frac{1}{4}$ & 0 & $\frac{3}{4}$ &  &  &    
    \end{sark_tableau}
    \label{tab:sark3_heun}
  }
  \caption{Three-stage SARK methods}
  \label{tab:sark3_examples}
\end{table}

The third order SARK methods in Table \ref{tab:sark3_examples} are
based on known third order Runge-Kutta methods with three stages. The additional coefficients $d_{ij}$ are chosen, such that (\ref{eq:cond_2})-(\ref{eq:cond_3}) are satisfied.

\subsection{Application of multiple steps within a tent}
\label{ssec:sark-mult-steps}

Recall that the ODE system we need to solve within one mapped tent
is~\eqref{eq:exactflow} for $0 < \hat t < 1.$ Since the $\hat t$
interval is not small, we subdivide it into $r$ subintervals and use
the previously described $s$-stage SARK scheme within each
subinterval, as described next. 

We subdivide the unit interval $[0,1]$ into $r$ subintervals
\[
  [\hat t_k, \hat t_{k+1}],
  \quad k=0,1, \ldots, r-1,
  \quad\text{ where } \hat t_k = \frac k r, 
\]
and apply~\eqref{eq:rktype} within each subinterval as described next.

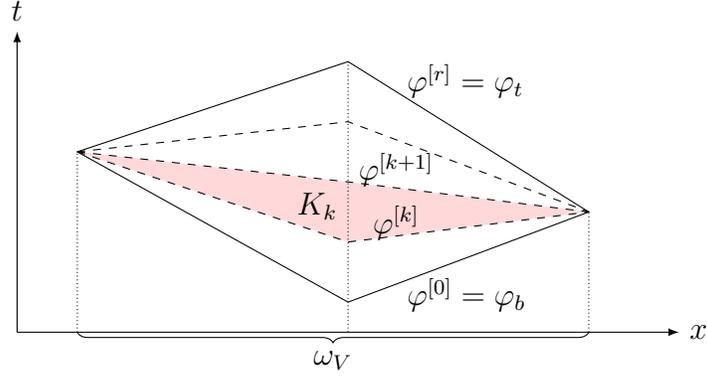
\begin{figure}
  \centering
  \begin{tikzpicture}[scale=8]
    \tikzstyle{redfill} = [fill=red!30,fill opacity=0.5]
    
    \coordinate (o) at (0,0);
    \coordinate (xshift) at (-0.05,0);    
    \coordinate (b1) at ($(0.15,0)+(xshift)$);
    \coordinate (b2) at ($(0.6,0)+(xshift)$);
    \coordinate (b3) at ($(1,0)+(xshift)$);
    
    \coordinate (yshift) at (0,-0.05);    
    \coordinate (tl) at ($(b1)+(yshift)+(0,0.35)$);
    \coordinate (tb) at ($(b2)+(yshift)+(0,0.1)$);
    \coordinate (tt) at ($(b2)+(yshift)+(0,0.5)$);
    \coordinate (tr) at ($(b3)+(yshift)+(0,0.25)$);
    
    \draw (tb) -- (tl) -- (tt) -- (tr) -- cycle;
    \node (phib) at ($0.8*(tb)+0.2*(tr)$)[right,yshift=-4pt] {$\vphi^{[0]}=\vphi_b$};
    \node (phit) at ($0.8*(tt)+0.2*(tr)$)[right,yshift=4pt] {$\vphi^{[r]}=\vphi_t$};
    \fill[redfill] ($0.75*(tb)+0.25*(tt)$) -- (tl) -- ($0.5*(tb)+0.5*(tt)$) -- (tr) -- cycle;
    \coordinate (phi1) at ($0.75*(tb)+0.25*(tt)$);
    \draw[dashed] (tl) -- (phi1) -- (tr);
    \node (phik) at ($0.8*(phi1)+0.2*(tr)$)[yshift=5pt] {$\vphi^{[k]}$};
    \coordinate (phi2) at ($0.5*(tb)+0.5*(tt)$);
    \draw[dashed] (tl) -- (phi2) -- (tr);
    \node (phik) at ($0.8*(phi2)+0.2*(tr)$)[yshift=7pt] {$\vphi^{[k+1]}$};
    \coordinate (phi3) at ($0.25*(tb)+0.75*(tt)$);
    \draw[dashed] (tl) -- (phi3) -- (tr);    
    \node at ($0.7*(tb)+0.3*(tt)$)[anchor=south east] {$K_k$};
    \draw[densely dotted] (b1) -- (tl);
    \draw[densely dotted] (b2) -- (tt);
    \draw[densely dotted] (b3) -- (tr);
    
    \draw [decorate,decoration={brace,amplitude=4pt,mirror}] (b1) -- (b3) node[midway,below,yshift=-3pt] {$\vp$};
        
    \draw[-latex] (o)--($(o)+(1.1,0)$) node[right] {$x$};
    \draw[-latex] (o)--($(o)+(0,0.5)$) node[above] {$t$};

  \end{tikzpicture}
  \caption{Illustration of the subtent $K_k$ (shaded) defined
    in~\eqref{eq:Kk}. It is the image under $\varPhi$ of the
    tensor product domain $\hat K_k = \vp \times (\that_k,\that_{k+1})$. \label{fig:map_tent}}
\end{figure}

First observe that the above splitting of the unit $\hat t$-interval
corresponds to subdividing the original tent $K$, as given
by~\eqref{eq:def_tent}, 
into $r$
``subtents'' (see Fig. \ref{fig:map_tent}) of the form
\begin{equation}
  \label{eq:Kk}
  K_k = \{ (x, t): \; x \in \vp, \;\varphi^{[k]} \le t \le \varphi^{[k+1]}\}
\end{equation}
where $ \varphi^{[k]} = \varphi(\hat t_k)$. Clearly $\varphi^{[0]}
= \varphi_b$ and $\varphi^{[r]} = \varphi_t$.

We then apply~\eqref{eq:rktype} to each of these
subtents. Accordingly, let $M_{0, [k]}$ be defined by~\eqref{eq:Mo}
after replacing $\varphi_b$ by $\varphi^{[k]}$. Keeping the
same definition of $A$ and $M_1$, let
$\tM_1^{[k]} = M_1 \circ M_{0, [k]}^{-1}$,
$\tA^{[k]} = A \circ M_{0, [k]}^{-1}$, and
$\tau^{[k]} = \hat t_{k+1} - \hat t_k$. Then the application
of~\eqref{eq:rktype} on each interval $[\hat t_k, \hat t_{k+1}]$
results in the following algorithm.

\begin{algorithm}\label{alg:SARKrs} \hfill
  \begin{enumerate}
  \item If the input is $Y_0$, an approximation to $Y(0)$ at the tent
    bottom, then set $Y^{[0]} = Y_0$.  If the input is $U_0$, an
    approximation to $U(0)$ at the tent bottom, then set
    $Y^{[0]} = Y_0 = M_0(U_0)$.
    
  \item   For $k = 0, 1, \ldots, r-1$ do:
    \begin{enumerate}
      
    \item For $i = 1, 2, \ldots, s$, compute 
      \[
        Z_i^{[k]}  = Y^{[k]} + \tau^{[k]} \sum_{j=1}^{i-1} d_{ij}\tM_1^{[k]}(Z_j^{[k]})
        + \tau^{[k]}
        \sum_{j=1}^{i-1} a_{ij}\tA^{[k]}(Z_j^{[k]}).
      \]
    \item  Compute 
      \[
        Y^{[k+1]}  = Y^{[k]} + \tau^{[k]} \sum_{i=1}^s b_i \tA^{[k]}(Z_i^{[k]}).
      \]
    \end{enumerate}
  \item Set
    \[
      Y_1^{r, s} = Y^{[r]}.
    \]
    Output this as the  approximation to $Y(1)$ at the tent top.
  \end{enumerate}
\end{algorithm}

We conclude this section by 
defining the propagation operators of the above
algorithm, which we shall use later. 
At step $k$, we define the (generally nonlinear)
partial propagation operator
$ T^{[k+1]} : V_h(\vpatch) \to V_h(\vpatch)$, 
using the intermediate quantities in the algorithm:
\begin{subequations}
  \label{eq:propagation-op}
  \begin{equation}
    \label{eq:prop_op_k}
      T^{[k+1]} (Y^{[k]}) =  Y^{[k+1]}.
  \end{equation}
Let the total propagation operator on the tent
$T : V_h(\vpatch) \to V_h(\vpatch)$ be defined by
\begin{equation}
  \label{eq:prop_op_tent}
  T = T^{[r]} \circ \ldots \circ T^{[2]} \circ T^{[1]}.
\end{equation}
\end{subequations}
Clearly, the input and output of the algorithm are related to $T$ by 
\begin{equation}
  \label{eq:prop_yrs}
  Y_1^{r, s}  = T(Y_0). 
\end{equation}

\section{Investigation of discrete stability} \label{sec:discrete-stability}

This section is devoted to remarks on the stability of the new SARK
schemes. While it is common to study stability of ODE solvers by
applying them to a simple scalar ODE, keeping our application of
spatially varying hyperbolic solutions in mind, we consider changes in
an energy-like measure on the solution $U(\hat t)$.  Recall that
$U(\that)\in\RRR^m$ is the coefficient vector of the basis expansion
of the mapped finite element solution $\hat u(x,\that)\in V_h(\vp)$,
as defined by~\eqref{eq:basis_exp}. We limit ourselves to the case where the
energy-like quantity
\begin{align}
  \label{eq:mnorm}
  \Mnorm{U(\that )}{\that}^2:= \int_\vp M(\hat t, U) \cdot \hat u =
  \int_\vp \left(M_0 (U) - \that M_1(U)\right)\cdot \hat u
\end{align}
is a {\em norm} and (the generally nonlinear operators) $M, M_0$ and $M_1$
defined in (\ref{eq:Mt}), (\ref{eq:Mo}) and (\ref{eq:M1}),
respectively, are linear, so that we may rewrite
$M(\hat t, U) = M(\hat t) U$ using the linear operator
$M(\hat t):= M_0 - \hat t M_1 :
\RRR^m \to V_h(\vpatch)$.
For many standard linear hyperbolic
systems, the causality condition can be used to easily show that
$M(\hat t)$ is identifiable with a symmetric positive definite
matrix so that 
\eqref{eq:mnorm} indeed defines a norm. In the special case of
$g(v) = v$, we note that on flat advancing fronts, where
$\vphi(x,\that)$ is independent of $x$ for some fixed $\that$,
\eqref{eq:mnorm} reduces to
\begin{align*}
  \Mnorm{U(\hat t)}{\that}^2=\int_\vp \hat u \cdot \hat u,
\end{align*}
so $\Mnorm{U(\hat t) }{\that}$ becomes the familiar spatial $L^2$ norm of $ \hat
u(\cdot, \hat t)$.

\subsection{Our procedure to study linear stability}

Stability of the scheme within a tent can be understood by studying
the discrete analogue of the ratio
$\Mnorm{U(1)}{1} / \Mnorm{U(0)}{0}$ for all possible initial data
$U(0)$. This amounts to studying the norm of the discrete propagation
operator for $U$, which we proceed to formulate.  First, recall the
connection between $U$ and $Y$, namely
$Y(\hat t) = M(\hat t) U(\hat t)$. Algorithm~\ref{alg:SARKrs}, takes
as input an approximation $U_0$ to $U(0)$ at the tent bottom and
outputs $Y_1^{r, s}$, an approximation to $Y(1)$ at the tent
top. Hence the associated approximation to $U(1)$ is
\[
  U^{r,s}_1 := M(1)^{-1} Y_1^{r, s}.
\]
Next, recall the discrete propagation operator defined
by~\eqref{eq:propagation-op}.  It is now a linear operator that maps
$Y_0= M(0) U_0$ to $Y_1^{r, s}$ according to~\eqref{eq:prop_yrs}.  Define the {\em
 tent propagation matrix} $S : \RRR^m \to \RRR^m$ by
\begin{equation}
  \label{eq:S}
  S = M(1)^{-1} T M(0).
\end{equation}
Clearly, \eqref{eq:prop_yrs} implies that
\begin{equation}
  \label{eq:prop_urs}
  U_1^{r,s} = S \,U_0.
\end{equation}

The discrete analogue of 
$\Mnorm{U(1)}{1} / \Mnorm{U(0)}{0}$ is 
$\Mnorm{U_1^{r, s}}{1} / \Mnorm{U_0}{0}$ which can be bounded using
the following   norm of $S$:
\[
  \| S \|_{L(M(0), M(1))} \,=
  \sup_{0 \ne W \in \RRR^m}
  \frac{ \| S W \|_{M(1)}} { \| W \|_{M(0)}}.
\]
It is immediate from~\eqref{eq:prop_urs} that
$\| U_1^{r,s}\|_{M(1)} \le \| S \|_{L(M(0), M(1))} \|
U_0\|_{M(0)}$. Thus the study of stability of SARK schemes is reduced
to computing estimates for the norm of $S$.

We now describe how we computed the norm of $S$ for some examples
below. Writing
$\hat v (x, \that ) = \sum_{i=1}^m V_i(\that) \psi_i(x)$ and
$\hat w (x, \that ) = \sum_{i=1}^m W_i(\that) \psi_i(x)$, in analogy
with the basis expansion of $\hat u$ in~\eqref{eq:basis_exp}, let
$\M_{\hat t}$ be the $m\times m$ symmetric positive definite matrix
satisfying
$W^\top \M_{\hat t} V = \int_\vpatch M(\hat t) V \cdot \hat w$. Then
\begin{align*}
  \| S \|^2_{L(M(0), M(1))}
  & = 
    \sup_{0 \ne W \in \RRR^m}
    \frac{ (S W)^\top \M_1 (SW)} { W^\top \M_0 W}
  \\
  & = 
    \sup_{0 \ne W \in \RRR^m}
    \frac{ W^\top (S^\top \M_1 S) W} { W^\top \M_0 W}
      \\
  & = 
    \sup\{ |\lambda|: \;\exists 0 \ne X \in \RRR^m \text{ satisfying }
    (S^\top \M_1 S) X  = \lambda \M_0 X\}.
\end{align*}
Thus, to investigate the stability of a scheme, we computed $T^{[k]}$
from
the scheme's Butcher-like tableau,
then $T$ by~\eqref{eq:prop_op_tent}, followed 
by $S$ per~\eqref{eq:S}, and
finally,  the square root of the spectral radius of
$ \M_0^{-1} (S^\top \M_1 S)$, which equals $\| S \|_{L(M(0), M(1))}$
as shown above.  We expand on the first of these steps in the next
few subsections by displaying $T^{[k]}$ for some SARK
schemes and end this section by reporting
our numerical estimates for $\| S \|_{L(M(0), M(1))}$  for an example.

\subsection{Propagation operator of two-stage SARK methods}
\label{ssec:propmat-sark2}
For an arbitrary two-stage SARK method the only non-zero coefficients
are $b_1, b_2, a_{21}, d_{21}$. For a given $Y^{[k]}=M_{0,[k]}U^{[k]}$
we obtain
\begin{align*}
  Z_1^{[k]} & = Y^{[k]}, \\  
  Z_2^{[k]} & = Y^{[k]} + \tau^{[k]} d_{21}\tM_1^{[k]} Z_1^{[k]}
              + \tau^{[k]} a_{21}\tA^{[k]} Z_1^{[k]}\\
            & = \left(I + \tau^{[k]} \left(d_{21}\tM_1^{[k]} + a_{21}\tA^{[k]}\right)\right) Y^{[k]},
\end{align*}
with the identity matrix $I\in\RRR^{m\times m}$.
The propagation from $\that_k$ to $\that_{k+1}$ reads
\begin{align*}
  Y^{[k+1]} & = Y^{[k]} + \tau^{[k]} \left(b_1 \tA^{[k]} Z_1^{[k]} + b_2 \tA^{[k]} Z_2^{[k]}\right)\\
            & = \left(I + \tau^{[k]}(b_1+b_2) \tA^{[k]}
              + (\tau^{[k]})^2 \tA^{[k]}\left(b_2d_{21}\tM_1^{[k]}
              + b_2a_{21} \tA^{[k]}\right)\right) Y^{[k]} \\
            & = \left(I + \tau^{[k]}\tA^{[k]}
              + \tfrac{1}{2}(\tau^{[k]})^2 \tA^{[k]}\left(\tM_1^{[k]} + \tA^{[k]} \right)\right) Y^{[k]},
\end{align*}
where we used the order conditions (\ref{eq:cond_1}) and
(\ref{eq:cond_2}) for second order methods. This results in the
propagation matrix
\begin{align*}
  T^{[k]} &= I + \tau^{[k]}\tA^{[k]}
                + \tfrac{1}{2}(\tau^{[k]})^2 \tA^{[k]}\left(\tM_1^{[k]} + \tA^{[k]} \right),
\end{align*}
such that $Y^{[k+1]} = T^{[k]}Y^{[k]}$.

\subsection{Propagation operator  of three-stage SARK methods}
\label{ssec:propmat-sark3}

A similar calculation for three-stage SARK methods, using the order
conditions (\ref{eq:cond_1})-(\ref{eq:cond_3}), leads to the
propagation matrix
\begin{align*}
  T^{[k]} =I &+ \tau^{[k]}\tA^{[k]}
                   + \tfrac{1}{2}(\tau^{[k]})^2 \tA^{[k]}\left(\tM_1^{[k]} + \tA^{[k]} \right)\\
                 & + \tfrac{1}{6}(\tau^{[k]})^3\tA^{[k]}\left(2\tM_1^{[k]} + \tA^{[k]} \right)
                   \left(\tM_1^{[k]} + \tA^{[k]} \right).
\end{align*}

\subsection{Discrete stability measure for a model problem}
\label{sec:discr-stab-conv}

We report the practically observed values of the previously described
stability measure (namely the norm $\| S \|_{L(M(0), M(1))}$)
for some SARK schemes
applied to the two-dimensional
convection equation
\begin{align*}
  \pd_t u(x,t) + \divx \left(b\,u(x,t)\right) = 0,\quad\forall(x,t)\in\Omega_0\times(0,t_{\max}],
\end{align*}
with $\Omega_0=[0,1]^2,~t_{\max}=0.05$, the flux field $b=(1,1)^\top$
and periodic boundary conditions. The time slab
$\Omega = \Omega_0\times(0,t_{\max})$ is filled with tents. Within
each such tent $K_i$, let $C_i$ denote the norm
$\| S \|_{L(M(0), M(1))}$ computed with $S, M(0),$ and $M(1)$ specific
to that tent. We expect $C_i$ to be close to one for a stable method.
Let   
\begin{align}  
  \label{eq:stab_cbar}
  \bar C := \max_{i} \; \left\{ C_i - 1 \right\}, 
\end{align}
where the maximum is taken over all tents in the time slab.
To gain an understanding of practical stability,
we examine the values of $\bar C$ as a function of the number of SARK stages
($s$), polynomial degree ($p$),  and more importantly,
the number of substeps per tent ($r$).

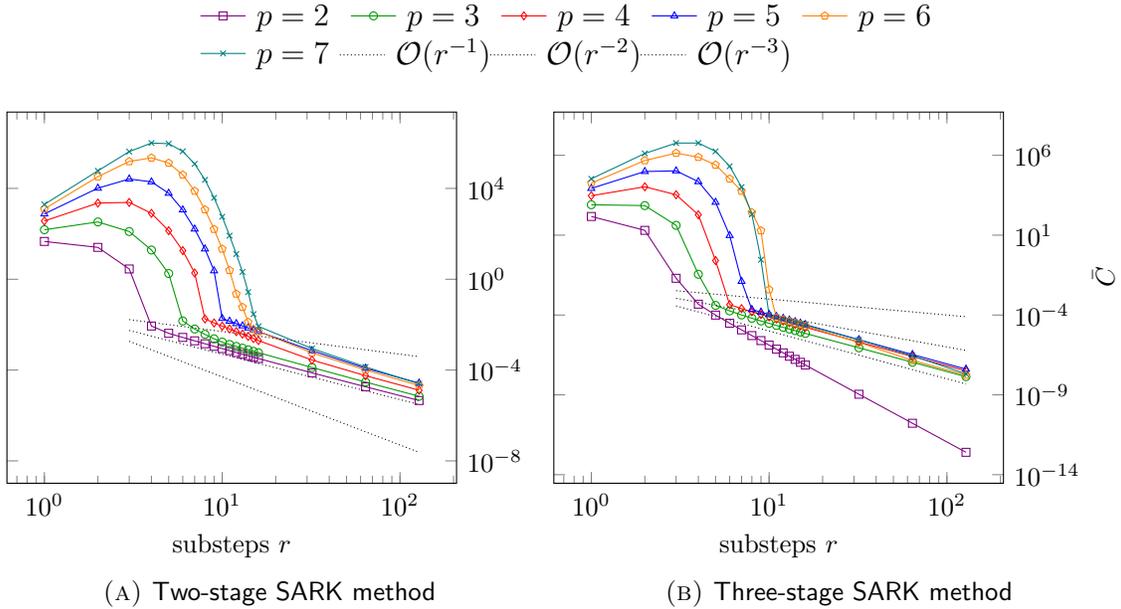
\begin{figure}[h]
  \captionsetup{width=.45\linewidth}
  \centering
  \begin{tikzpicture}
    \tikzset{mark options={mark size=1.5}}%
    \def\xlabeloffset{2};
    \def\ylabeloffset{-0.5};
    \node (a) at (0,0) {\ref{pgfplots:sark2_p2} $p=2$};
    \node at ($(a)+(\xlabeloffset,0)$) {\ref{pgfplots:sark2_p3} $p=3$};
    \node at ($(a)+(2*\xlabeloffset,0)$) {\ref{pgfplots:sark2_p4} $p=4$};
    \node at ($(a)+(3*\xlabeloffset,0)$) {\ref{pgfplots:sark2_p5} $p=5$};
    \node at ($(a)+(4*\xlabeloffset,0)$) {\ref{pgfplots:sark2_p6} $p=6$};
    \node at ($(a)+(0,\ylabeloffset)$) {\ref{pgfplots:sark2_p7} $p=7$};
    \node at ($(a)+(\xlabeloffset,\ylabeloffset)$) {\ref{pgfplots:Or} $\cO(r^{-1})$};
    \node at ($(a)+(2*\xlabeloffset,\ylabeloffset)$) {\ref{pgfplots:Or2} $\cO(r^{-2})$};
    \node at ($(a)+(3*\xlabeloffset,\ylabeloffset)$) {\ref{pgfplots:Or3} $\cO(r^{-3})$};
  \end{tikzpicture}
  \pgfplotscreateplotcyclelist{sarkplotstyles}{%
    {color=violet,mark=square},
    {green!60!black,mark=o},
    {red,mark=diamond},
    {color=blue,mark=triangle},
    {color=orange,mark=pentagon},
    {color=teal,mark=x}}
  \subfloat[Two-stage SARK method]
  {
    \label{fig:conv2d_stab_sark2}
    \begin{tikzpicture}
      \tikzset{mark options={mark size=1.5}}%
      \begin{loglogaxis}[%
        width=0.48\textwidth,%
        tick label style={font=\footnotesize},%
        yticklabel pos=right,
        xlabel={\footnotesize substeps $r$},%
        cycle list name = sarkplotstyles,%
        ]%
        \addplot table[x=r,y=lammax_p2_ws4]{table_conv2d_stability_sark2_ws4_fine.txt};
        \label{pgfplots:sark2_p2}
        \addplot table[x=r,y=lammax_p3_ws4]{table_conv2d_stability_sark2_ws4_fine.txt};
        \label{pgfplots:sark2_p3}
        \addplot table[x=r,y=lammax_p4_ws4]{table_conv2d_stability_sark2_ws4_fine.txt};
        \label{pgfplots:sark2_p4}
        \addplot table[x=r,y=lammax_p5_ws4]{table_conv2d_stability_sark2_ws4_fine.txt};
        \label{pgfplots:sark2_p5}
        \addplot table[x=r,y=lammax_p6_ws4]{table_conv2d_stability_sark2_ws4_fine.txt};
        \label{pgfplots:sark2_p6}
        \addplot table[x=r,y=lammax_p7_ws4]{table_conv2d_stability_sark2_ws4_fine.txt};
        \label{pgfplots:sark2_p7}
        
        \addplot[densely dotted] table[restrict expr to domain={\coordindex}{2:18},x=r,y expr={0.05/x}] {table_conv2d_stability_sark2_ws4_fine.txt};
        \label{pgfplots:Or}
        \addplot[densely dotted] table[restrict expr to domain={\coordindex}{2:18},x=r,y expr={0.05/x^2}] {table_conv2d_stability_sark2_ws4_fine.txt};
        \label{pgfplots:Or2}
        \addplot[densely dotted] table[restrict expr to domain={\coordindex}{2:18},x=r,y expr={0.05/x^3}] {table_conv2d_stability_sark2_ws4_fine.txt};
        \label{pgfplots:Or3}
      \end{loglogaxis}
    \end{tikzpicture}
  }
  \subfloat[Three-stage SARK method]
  {
    \label{fig:conv2d_stab_sark3}
    \begin{tikzpicture}
      \tikzset{mark options={mark size=1.5}}%
      \begin{loglogaxis}[%
        width=0.48\textwidth,%
        tick label style={font=\footnotesize},%
        yticklabel pos=right,
        xlabel={\footnotesize substeps $r$},%
        ylabel={\footnotesize $\bar C$},%
        ylabel style = {at={(1.38,0.5)} ,anchor=north west},%
        cycle list name = sarkplotstyles,%
        ]%
        \addplot table[x=r,y=lammax_p2_ws4]{table_conv2d_stability_sark3_ws4_fine.txt};
        \addplot table[x=r,y=lammax_p3_ws4]{table_conv2d_stability_sark3_ws4_fine.txt};
        \addplot table[x=r,y=lammax_p4_ws4]{table_conv2d_stability_sark3_ws4_fine.txt};
        \addplot table[x=r,y=lammax_p5_ws4]{table_conv2d_stability_sark3_ws4_fine.txt};
        \addplot table[x=r,y=lammax_p6_ws4]{table_conv2d_stability_sark3_ws4_fine.txt};
        \addplot table[x=r,y=lammax_p7_ws4]{table_conv2d_stability_sark3_ws4_fine.txt};
        
        \addplot[densely dotted] table[restrict expr to domain={\coordindex}{2:18},x=r,y expr={0.01/x}] {table_conv2d_stability_sark3_ws4_fine.txt};
        \addplot[densely dotted] table[restrict expr to domain={\coordindex}{2:18},x=r,y expr={0.01/x^2}] {table_conv2d_stability_sark3_ws4_fine.txt};
        \addplot[densely dotted] table[restrict expr to domain={\coordindex}{2:18},x=r,y expr={0.01/x^3}] {table_conv2d_stability_sark3_ws4_fine.txt};
      \end{loglogaxis}
    \end{tikzpicture}
  }
  \caption{Observed dependence of $\bar C$ on $r$ for $p=2,3,4,5,6,7$ and
    $s=2, 3$.}
  \label{fig:conv2d_stab}
\end{figure}

In all our numerical experiments, we observed that on each tent, for a
fixed $s$, the norm $\| S \|_{L(M(0), M(1))}$ tends to $1$ with
increasing number of substeps $r$, and moreover, we discovered a
dependence of the following form
\begin{align*}
  \| S \|_{L(M(0), M(1))} = 1 + \mathcal{O}\left(r^{-s}\right)
\end{align*}
on each tent $K_i$. Therefore, we organize our report on numerical
stability observations into 
plots of   values of $\bar C$ as a function of $r$. We do so
for two SARK methods, one with $s=2$
and another with $s=3$. The results are displayed  in Fig.~\ref{fig:conv2d_stab}.
After a prominent preasymptotic
region, we observe that $\bar C$,
as a function of $r$, exhibits the rate
$\mathcal{O}\left(r^{-s}\right)$ in all cases, except  one.

The exceptional case is the case $p=2$ in
Fig.~\ref{fig:conv2d_stab_sark3}, where the stability measure
approaches the ideal value of $1$ much faster. We do not have an
explanation for this observation.

Note that all the plotted curves in Fig.~\ref{fig:conv2d_stab} shift
to the top and right as $p$ increases, i.e., the 
number of substeps~$r$ required to keep the same stability measure 
$\bar C$ increases with $p$. This behavior is akin to the 
$p$-dependence of the CFL-conditions of standard time stepping schemes.

\section{Numerical results}  \label{sec:numerical}

In this section, we collect our observations on
the performance of the new SARK
schemes, on the one-dimensional Burger's equation (in
\S\ref{ssec:conv-rates-burg}) and the two-dimensional Euler system
(in
\S\ref{ssec:convergence-rates-euler}--\S\ref{ssec:mach-3-forward}).
While \S\ref{ssec:convergence-rates-euler}
focuses on the
study of convergence rates for a smooth Euler solution,
\S\ref{ssec:mach-3-forward} presents the application of SARK scheme on
the computationally challenging problem of simulating a 
Mach~3  wind tunnel with a forward-facing step.

\subsection{Convergence rates for Burger's equation}
\label{ssec:conv-rates-burg}

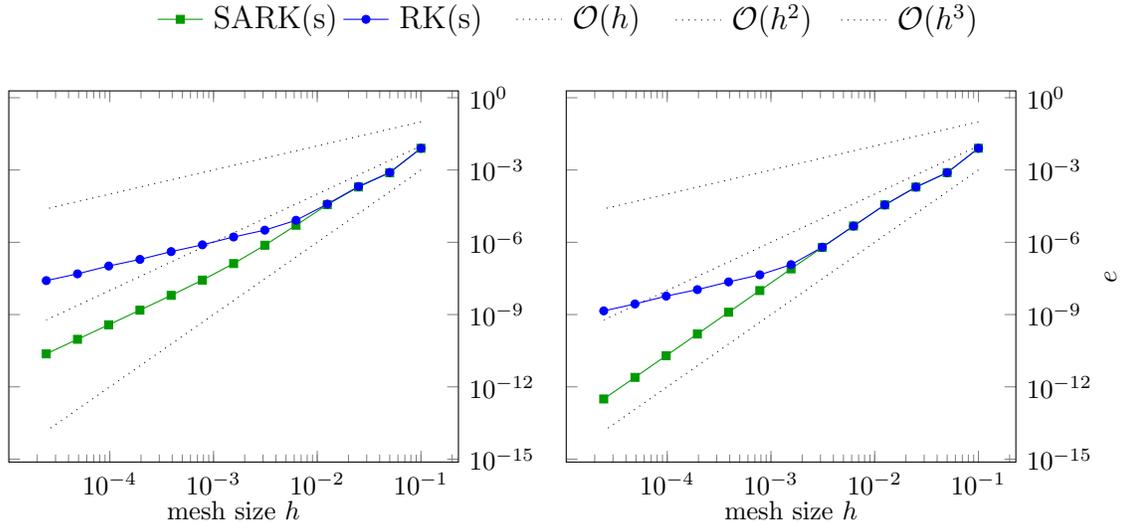
\begin{figure}[b] 
  \captionsetup{width=.45\linewidth}
  \centering
  \begin{tikzpicture}
    \tikzset{mark options={mark size=1.5}}%
    \def\labeloffset{2.2};
    \node (a) at (0,0) {\ref{pgfplots:sark} SARK(s)};
    \node at ($(a)+(\labeloffset,0)$) {\ref{pgfplots:rk} RK(s)};
    \node at ($(a)+(2*\labeloffset,0)$) {\ref{pgfplots:Oh} $\cO(h)$};
    \node at ($(a)+(3*\labeloffset,0)$) {\ref{pgfplots:Oh2} $\cO(h^2)$};
    \node at ($(a)+(4*\labeloffset,0)$) {\ref{pgfplots:Oh3} $\cO(h^3)$};
  \end{tikzpicture}
  \subfloat[Convergence rates obtained from  SARK(2, Ralston) method
  (see Table \ref{tab:sark2_ralston}) and the standard  Ralston method.]
  {
    \label{fig:burgers1d_sark2_rates}
    \begin{tikzpicture}
      \tikzset{mark options={mark size=1.5}}%
      \begin{loglogaxis}[%
        width=0.48\textwidth,%
        tick label style={font=\footnotesize},%
        yticklabel pos=right,
        xlabel={\footnotesize mesh size $h$},%
        xlabel style = {at={(0.5,0.03)} ,anchor=north},%
        ]%
        \addplot[color=green!60!black,mark=square*] table[x=h,y=l2e_SARK2_ralston] {table_burgers1d_conv_sark.txt};
        \label{pgfplots:sark}
        \addplot[color=blue,mark=*] table[x=h,y=l2e_RK2_ralston] {table_burgers1d_conv_rk.txt};
        \label{pgfplots:rk}

        \addplot[dotted] table[x=h,y=h] {table_burgers1d_conv_sark.txt};
        \label{pgfplots:Oh}
        \addplot[dotted] table[x=h,y expr={x^2}] {table_burgers1d_conv_sark.txt};
        \label{pgfplots:Oh2}
        \addplot[dotted] table[x=h,y expr={x^3}] {table_burgers1d_conv_sark.txt};
        \label{pgfplots:Oh3}
      \end{loglogaxis}
    \end{tikzpicture}
  }
  \subfloat[Convergence rates obtained from  SARK(3, Heun) method (see Table
  \ref{tab:sark3_heun}) and the standard Heun scheme.]
  {
    \label{fig:burgers1d_sark3_rates}
    \begin{tikzpicture}
      \tikzset{mark options={mark size=1.5}}%
      \begin{loglogaxis}[%
        width=0.48\textwidth,%
        tick label style={font=\footnotesize},%
        yticklabel pos=right,
        xlabel={\footnotesize mesh size $h$},%
        xlabel style = {at={(0.5,0.03)} ,anchor=north},%
        ylabel={\footnotesize $e$},%
        ylabel style = {at={(1.38,0.5)} ,anchor=north},%
        ]%
        \addplot[color=green!60!black,mark=square*] table[x=h,y=l2e_SARK3_heun] {table_burgers1d_conv_sark.txt};
        \addplot[color=blue,mark=*] table[x=h,y=l2e_RK3_heun] {table_burgers1d_conv_rk.txt};

        \addplot[dotted] table[x=h,y=h] {table_burgers1d_conv_sark.txt};
        \addplot[dotted] table[x=h,y expr={x^2}] {table_burgers1d_conv_sark.txt};
        \addplot[dotted] table[x=h,y expr={x^3}] {table_burgers1d_conv_sark.txt};
      \end{loglogaxis}
    \end{tikzpicture}
  }
  \caption{Plots of the error $e$ defined in
    (\ref{eq:burgers1d_example}) for SARK and RK methods.}
  \label{fig:burgers1d_sark_rates}
\end{figure}

Let us begin by returning to the one-dimensional model problem of
Section~\ref{sec:prob} to show that the SARK methods do
{\em not} suffer from the previously described convergence order reduction.
For this discussion,  the equation and error $e$ are as in
(\ref{eq:burgers1d_example}). We apply Algorithm \ref{alg:SARKrs} with 
SARK schemes of $s=2$ and $s=3$ stages, collect values
of $e$ for various $h$ and plot them in
Fig.~\ref{fig:burgers1d_sark_rates}.

The data shown in Fig.~\ref{fig:burgers1d_sark_rates} was generated
with the polynomial order $p=2$ in space and $h=2^{-i} / 10$ for
$i=0\dots 12$. The 
tents were built so
that~(\ref{eq:causality_constraint}) is satisfied with
$c_{\max}=8$.
Algorithm \ref{alg:SARKrs} is applied with 
$r=4$ substeps within each tent. As $h$ decreases, in
Fig.~\ref{fig:burgers1d_sark2_rates} 
we eventually see quadratic
convergence for the two-stage SARK method (although the convergence
rate seems to be slightly higher in a preasymptotic regime), 
while the rate of the underlying
standard Runge-Kutta method drops to first order. The three-stage SARK method in
Fig.~\ref{fig:burgers1d_sark3_rates} shows cubic convergence while the
rate of the 
underlying standard Runge-Kutta method drops to first order again.
These plots clearly show the benefit of using SARK scheme over
the corresponding standard Runge-Kutta scheme. 

\subsection{Convergence rates for a 2D Euler system}
\label{ssec:convergence-rates-euler}

Now we apply SARK methods to the Euler
system. Similar to the Burger's example, which we discussed in the
previous section, we choose smooth initial data and fix a  final time
before the onset of shock so that no limiting is
needed. 

\begin{figure}
  \captionsetup{width=.45\linewidth}
  \centering
  \begin{tikzpicture}
    \tikzset{mark options={mark size=1.5}}%
    \def\labeloffset{2.2};
    \node (a) at (0,0) {\ref{pgfplots:euler2d_sark} SARK(s)};
    \node at ($(a)+(\labeloffset,0)$) {\ref{pgfplots:euler2d_rk} RK(s)};
    \node at ($(a)+(2*\labeloffset,0)$) {\ref{pgfplots:euler2d_Oh} $\cO(h)$};
    \node at ($(a)+(3*\labeloffset,0)$) {\ref{pgfplots:euler2d_Oh2} $\cO(h^2)$};
    \node at ($(a)+(4*\labeloffset,0)$) {\ref{pgfplots:euler2d_Oh3} $\cO(h^3)$};
  \end{tikzpicture}
  \subfloat[Convergence rates obtained from  SARK(2, Ralston) method
  (see Table \ref{tab:sark2_ralston}) and the standard  Ralston method.]
  {
    \label{fig:euler2d_sark2_rates}
    \begin{tikzpicture}
      \tikzset{mark options={mark size=1.5}}%
      \begin{loglogaxis}[%
        width=0.48\textwidth,%
        tick label style={font=\footnotesize},%
        yticklabel pos=right,
        xlabel={\footnotesize dof},%
        xlabel style = {at={(0.5,0.03)} ,anchor=north},%
        ]%
        \addplot[color=green!60!black,mark=square*] table[x=ndofs,y=l2error] {table_euler2d_conv_sark2ralston.txt};
        \label{pgfplots:euler2d_sark}
        \addplot[color=blue,mark=*] table[x=ndofs, y=l2error] {table_euler2d_conv_rk2ralston.txt};
        \label{pgfplots:euler2d_rk}

        \addplot[dotted] table[x=ndofs,y expr={0.01/x^(1/2)}] {table_euler2d_conv_sark2ralston.txt};
        \label{pgfplots:euler2d_Oh}
        \addplot[dotted] table[x=ndofs,y expr={0.2/x}] {table_euler2d_conv_sark2ralston.txt};
        \label{pgfplots:euler2d_Oh2}
        \addplot[dotted] table[x=ndofs,y expr={3/x^(3/2)}] {table_euler2d_conv_sark2ralston.txt};
        \label{pgfplots:euler2d_Oh3}
      \end{loglogaxis}
    \end{tikzpicture}
  }
  \subfloat[Convergence rates obtained from  SARK(3, Heun) method (see Table
  \ref{tab:sark3_heun}) and the standard Heun scheme.]
  {
    \label{fig:euler2d_sark3_rates}
    \begin{tikzpicture}
      \tikzset{mark options={mark size=1.5}}%
      \begin{loglogaxis}[%
        width=0.48\textwidth,%
        tick label style={font=\footnotesize},%
        yticklabel pos=right,
        xlabel={\footnotesize dof},%
        xlabel style = {at={(0.5,0.03)} ,anchor=north},%
        ylabel={\footnotesize $e$},%
        ylabel style = {at={(1.38,0.5)} ,anchor=north},%
        ]%
        \addplot[color=green!60!black,mark=square*] table[x=ndofs,y=l2error] {table_euler2d_conv_sark3heun.txt};
        \addplot[color=blue,mark=*] table[x=ndofs,y=l2error] {table_euler2d_conv_rk3heun.txt};
        
        \addplot[dotted] table[x=ndofs,y expr={0.1/x^(1/2)}] {table_euler2d_conv_sark3heun.txt};
        \addplot[dotted] table[x=ndofs,y expr={2.5/x}] {table_euler2d_conv_sark3heun.txt};
        \addplot[dotted] table[x=ndofs,y expr={50/x^(3/2)}] {table_euler2d_conv_sark3heun.txt};
      \end{loglogaxis}
    \end{tikzpicture}
  }
  \caption{Error $e$ as defined in (\ref{eq:euler2d_l2e}) over spatial
    degrees of freedom (dof) for SARK and standard RK methods applied
    to the Euler equation on tents as described in
    (\ref{eq:euler2d_example}).}
  \label{fig:euler2d_sark_rates}
\end{figure}
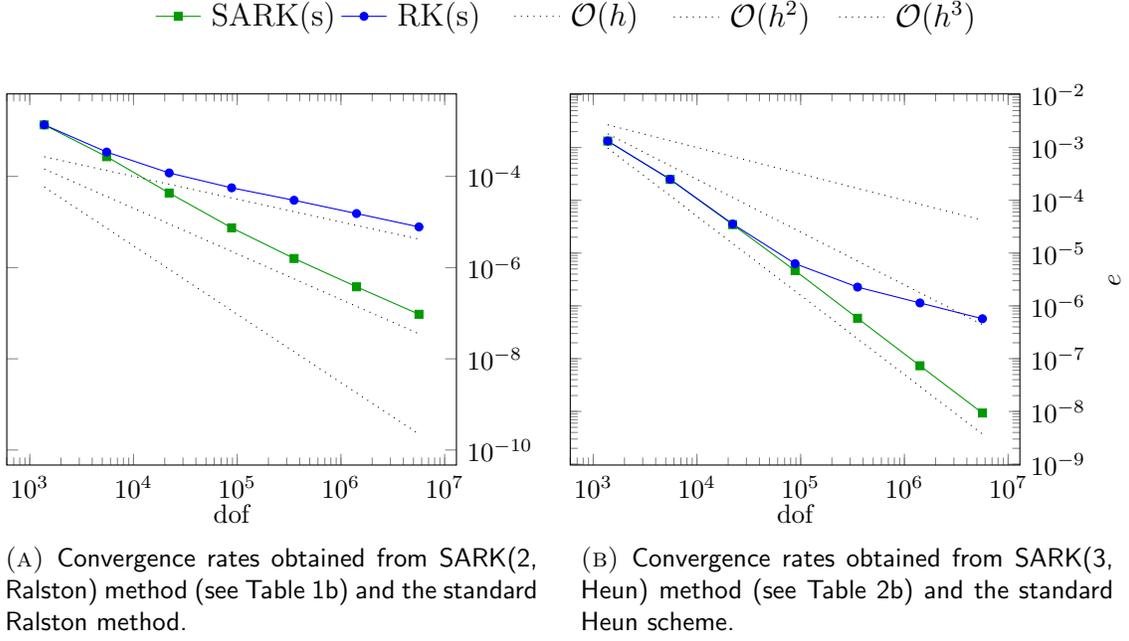

Recall that the Euler system fits into (\ref{eq:conslaw}) with 
\begin{subequations}
  \label{eq:euler2d_example}
  \begin{align}
    \label{eq:euler_system}
    u = \begin{pmatrix} \rho \\ m \\ E \end{pmatrix}, \quad g(u) = u, \quad
    f(u) = \begin{pmatrix} m \\ m\otimes m/\rho + \pp I \\ ( E + \pp)  m / \rho \end{pmatrix}.
  \end{align}
Here   
the functions $\rho:\Omega_0\to\RRR$,
  $m:\Omega_0\to\RRR^2$ and $E:\Omega_0\to\RRR$ denote the density,
  momentum, and total energy of a perfect gas in the spatial domain
  $\Omega_0=[0,1]^2$. Furthermore, we use $\pp=\tfrac{1}{2}\rho \tT$
  for the pressure,
  $\tT = \tfrac{4}{d}\big(
    \tfrac{E}{\rho}-\tfrac{1}{2}\tfrac{|m|^2}{\rho^2}\big)$ for the
  temperature and $d=5$ denotes the degrees of freedom of the gas
  particles. The initial values are set by
  \begin{align}
    \rho_0 &= 1+e^{-100((x-0.5)^2+(y-0.5)^2)}, \\
    m_0 &= (0,0)^\top, \\
    \pp_0 &= 1+e^{-100((x-0.5)^2+(y-0.5)^2)},
  \end{align}
  and the final time $t_{\max}=0.1$.

  The data shown in Fig. \ref{fig:euler2d_sark_rates} was generated
  with polynomial degree $p=2$ in space and mesh sizes
  $h=0.1\times 2^{-i}$, for $i=0\dots 6$. For the tent generation
  $c_{\max}$ in (\ref{eq:causality_constraint}) was set to $8$ and the
  number of substeps $r=4$. Since we do not have an exact solution in
  closed form, we compare 
  the numerical solution
  computed using $c_{\max}$ with a ``reference
  solution'' computed with the higher
  characteristic speed $2\cdot c_{\max}$. The latter requires many
  more tents to reach the final time. Let the former and latter
  approximations to $u(\cdot, t_{\max})$ be denoted by $u_h$ and
  $u_h^{\text{ref}}$, respectively.
  We define the error by
  \begin{align}
    \label{eq:euler2d_l2e}
    e:=\left\| u_h - u_h^{\text{ref}}\right\|_{L^2(\om_0)}.
  \end{align}  
\end{subequations}
This is the quantity that is plotted in Fig. \ref{fig:euler2d_sark_rates}.

The errors of the two-stage SARK method and the underlying RK method
is seen to diverge already for the first
refinement level
in Fig.~\ref{fig:euler2d_sark2_rates}.
While the SARK method shows the expected second
convergence order, the rate of the RK method drops to first order.
For the three-stage methods in Fig. \ref{fig:burgers1d_sark3_rates},
we see cubic convergence for both method for the first few refinements.
The convergence rate of the RK method eventually drops to first
order while the SARK converges at third order.

\subsection{Mach~3 wind tunnel}
\label{ssec:mach-3-forward}

We conclude with the well-known benchmark example \cite{WoodwColel84} of
the wind tunnel with a  forward-facing step onto which gas flows at 
Mach~3. The situation is modeled by the
already described Euler system  (\ref{eq:euler_system}), but now
with the initial values
\begin{align}
  \rho_0 = 1.4, \quad m_0 = \rho_0(3,0)^\top, \quad \pp_0 = 1
\end{align}
on a spatial domain $\Omega_0$ with a re-entrant corner at the
edge of the forward-facing step --- the domain and the boundary
conditions are exactly as illustrated in numerous previous works, see 
e.g., \cite[Fig.~4(a)]{mtp}.
Our numerical experience with this problem shows that it is
beneficial to use
high order local time stepping. As in our prior study~\cite{mtp},  we
use a spatially refined mesh near the re-entrant corner and let the
tents adapt, providing automatic local time stepping. In contrast
to the standard time stepping used in~\cite{mtp},
we now use one of the newly proposed  SARK schemes.

We shall
apply the SARK(3, Heun) method. Unlike the study in
\S\ref{ssec:convergence-rates-euler}, now we must 
handle multiple shocks that develop over time, so it is necessary to
add some 
stabilization  to the system. This is done by adding artificial
viscosity based on the entropy residual as suggested
by~\cite{GuermPasquPopov11}---details of this stabilization on tents
are exactly as already described  in~\cite{mtp}, so we omit them here.

\begin{figure}
  \centering
  \pgfplotsset{
    colormap={netgen32}{
      rgb255=(0, 0, 255),
      rgb255=(0, 32, 255),
      rgb255=(0, 65, 255),
      rgb255=(0, 98, 255),
      rgb255=(0, 131, 255),
      rgb255=(0, 164, 255),
      rgb255=(0, 197, 255),
      rgb255=(0, 230, 255),
      rgb255=(0, 255, 246),
      rgb255=(0, 255, 213),
      rgb255=(0, 255, 180),
      rgb255=(0, 255, 148),
      rgb255=(0, 255, 115),
      rgb255=(0, 255, 82),
      rgb255=(0, 255, 49),
      rgb255=(0, 255, 16),
      rgb255=(16, 255, 0),
      rgb255=(49, 255, 0),
      rgb255=(82, 255, 0),
      rgb255=(115, 255, 0),
      rgb255=(148, 255, 0),
      rgb255=(180, 255, 0),
      rgb255=(213, 255, 0),
      rgb255=(246, 255, 0),
      rgb255=(255, 230, 0),
      rgb255=(255, 197, 0),
      rgb255=(255, 164, 0),
      rgb255=(255, 131, 0),
      rgb255=(255, 98, 0),
      rgb255=(255, 65, 0),
      rgb255=(255, 32, 0),
      rgb255=(255, 0, 0),
    },
    colormap name = netgen32,
  }
  \subfloat[Logarithmic density]{
    \label{fig:mach3wt_density}
    \begin{tikzpicture}
      \node at (0,0.35) [anchor=south west]{%
        \includegraphics[width=0.55\paperwidth, trim=50 235 50 235, clip=true]{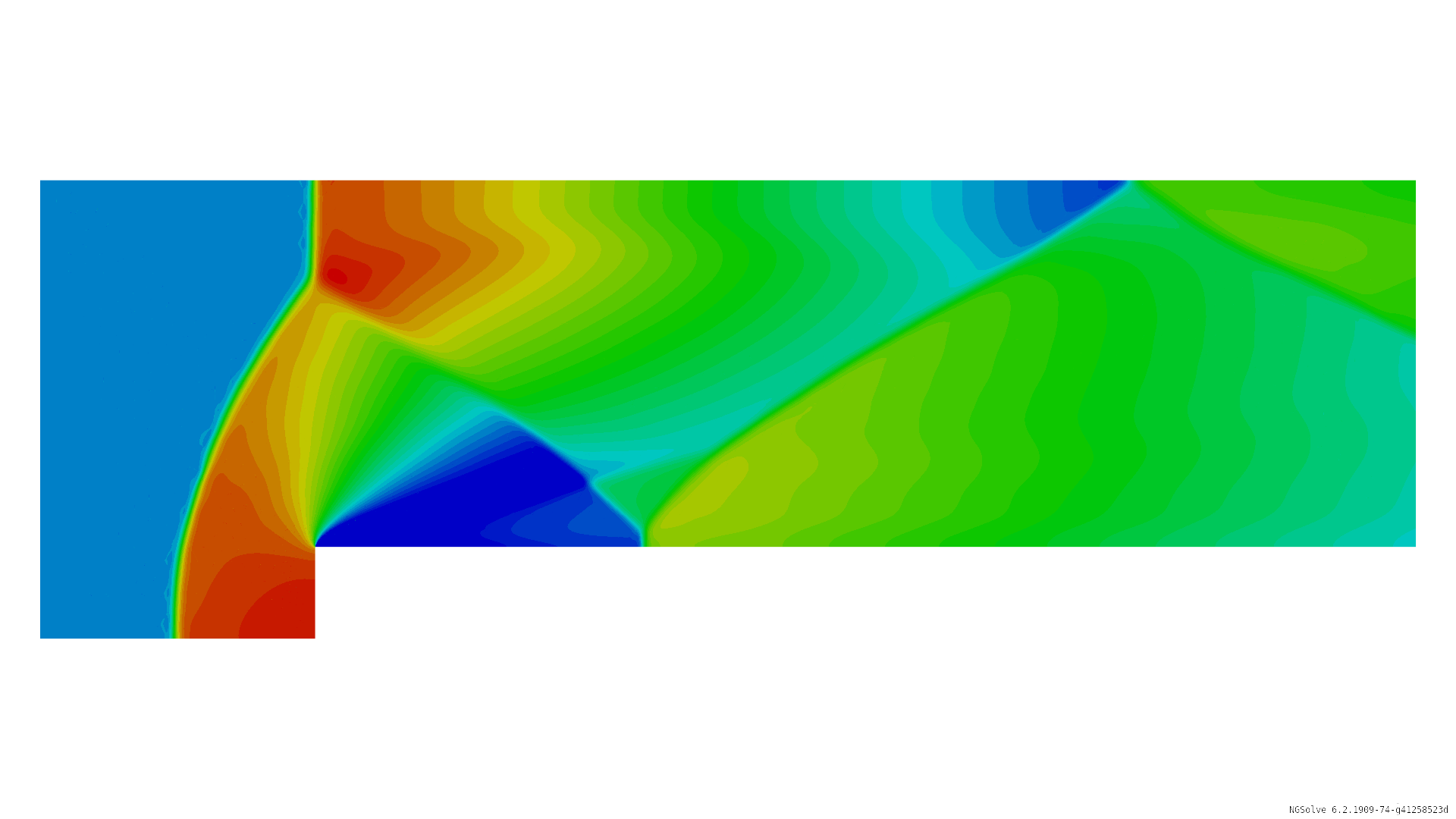}};
      \node at (3,0) [anchor=south west]{%
        \pgfplotscolorbardrawstandalone[
        colormap access=piecewise const,
        colorbar sampled,
        colorbar horizontal,
        point meta min=0,
        point meta max=1.85,
        colorbar style={
          axis line style={draw=none},
          width=0.25\paperwidth,
          height=10pt,
          xtick={0,0.45,0.9,1.35,1.8},
          xticklabel style = {font=\scriptsize},
        }]
        };
    \end{tikzpicture}
  }\\
  \subfloat[Entropy viscosity coefficient $\nu$]{
    \label{fig:mach3wt_visc}
    \begin{tikzpicture}
      \node at (0,0.35) [anchor=south west]{%
        \includegraphics[width=0.55\paperwidth, trim=50 235 50 235, clip=true]{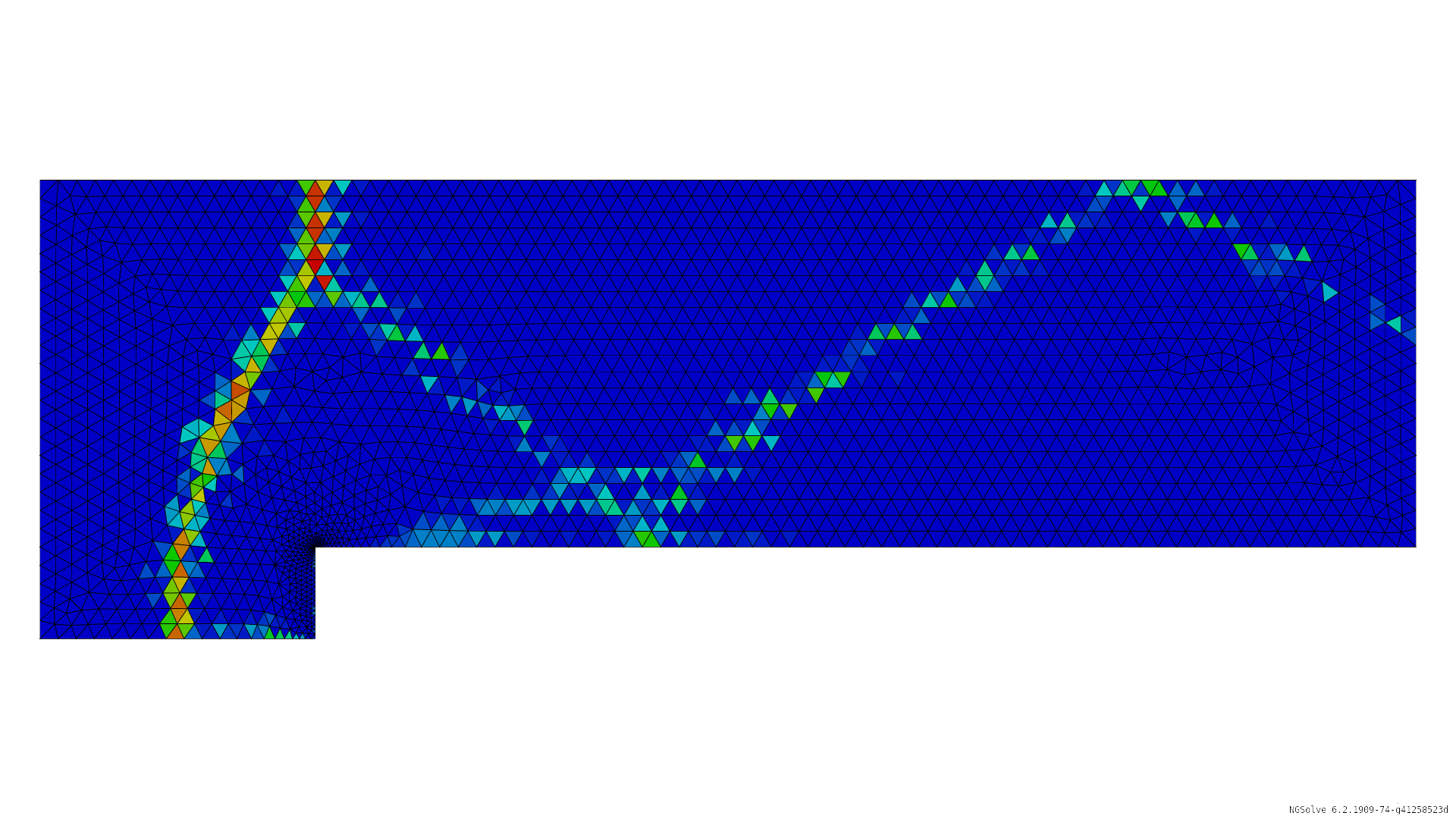}};
      \node at (3,0) [anchor=south west]{%
        \pgfplotscolorbardrawstandalone[
        colormap access=piecewise const,
        colorbar sampled,
        colorbar horizontal,
        point meta min=0,
        point meta max=0.008,
        colorbar style={
          axis line style={draw=none},
          width=0.25\paperwidth,
          height=10pt,
          xtick={0,0.002,0.004,0.006,0.008},
          xticklabel style = {font=\scriptsize},
          x tick scale label style={
            at={(xticklabel cs:1.1,-1em)},
            anchor=near xticklabel,
          },
        }]
        };
    \end{tikzpicture}
  }
  \caption{Solution of the Mach 3 wind tunnel with a forward-facing step at the final
    time $t_{\max}=4$ solved on 4128 triangles with SARK(3, Heun) and
    spatial degree $p=4$. }
  \label{fig:mach3wt}
\end{figure}

One of the components of the computed solution is
shown in Fig.~\ref{fig:mach3wt}. This  was generated with polynomial
order $p=4$ in space, maximal characteristic speed $c_{\max}=10$ and
$r=16$ substeps within each tent. The solution component (logarithmic
density) shown in Fig.~\ref{fig:mach3wt_density} is comparable with the
solution we previously obtained using standard methods in~\cite{mtp},
but now due to the higher accuracy of the new
SARK time
integration, we obtained a similar quality solution faster
(with the overall simulation time on the same processor
reduced by a factor of 10).
We also observed that the entropy residuals calculated off the computed
solution with SARK schemes led to a significantly 
reduced addition of  artificial viscosity. The artificial viscosity coefficients
generated by the entropy residual are shown in
Fig.~\ref{fig:mach3wt_visc}, which is about half the size of what is shown in
the corresponding plot in our earlier work~\cite[Fig.~5]{mtp}.

%
%

\bibliographystyle{spmpsci}      
\bibliography{waves_proc_references}   




\end{document}